\documentclass[11pt]{amsart}
\usepackage[linkcolor=blue, citecolor = blue]{hyperref}
\hypersetup{colorlinks=true}

\usepackage{amsthm,amsfonts,amsmath,amssymb,amscd} 

\usepackage{verbatim}
\usepackage{float}
\usepackage{wrapfig}
\usepackage{mathtools}
\usepackage[color=orange!20, linecolor=orange]{todonotes}

\usepackage{enumitem}

\usepackage{ytableau}

\usepackage{tikz}
\usetikzlibrary{decorations.markings}
\tikzstyle{vertex}=[circle, draw, inner sep=0pt, minimum size=4pt]

\usetikzlibrary{arrows,automata}
\usepackage{tikz, tikz-3dplot, pgfplots}
\usepackage{tkz-graph}
\usetikzlibrary[positioning,patterns]

\newtheorem{theorem}{Theorem}[section]
\newtheorem{proposition}[theorem]{Proposition}
\newtheorem{lemma}[theorem]{Lemma}
 
\theoremstyle{definition}
\newtheorem{corollary}[theorem]{Corollary}
\newtheorem{definition}[theorem]{Definition}

\theoremstyle{remark}
\newtheorem{remark}[theorem]{Remark}

\newcommand{\sgn}{\mathrm{sgn}}

\title[Grothendieck polynomials, Cauchy identities, and Dual filtered graphs]{Symmetric Grothendieck polynomials, skew Cauchy identities, and dual filtered Young graphs}

\author[Damir Yeliussizov]{Damir Yeliussizov}
\address{Department of Mathematics, UCLA, Los Angeles, CA 90095}
\email{\href{mailto:damir@math.ucla.edu}{damir@math.ucla.edu}}

\begin{document}

\begin{abstract}
Symmetric Grothendieck polynomials are analogues of Schur polynomials in the K-theory of Grassmannians. We build dual families of symmetric Grothendieck polynomials using Schur operators. With this approach we prove skew Cauchy identity and then derive various applications: skew Pieri rules, dual filtrations of Young's lattice, generating series and enumerative identities. We also give a new  explanation of the finite expansion property for products of Grothendieck polynomials. 
\end{abstract}

\maketitle

\section{Introduction}
Symmetric Grothendieck polynomials, also known as stable Grothendieck polynomials, are certain K-theoretic deformations of Schur functions. These functions were first studied by Fomin and Kirillov \cite{fk} as a stable limit of more general Grothendieck polynomials that generalize Schubert polynomials in another direction.

The symmetric Grothendieck polynomial $G_{\lambda}$ can be defined by the following combinatorial formula due to Buch \cite{buch} %as the following generating series
$$
G_{\lambda}(x_1, x_2, \ldots) = \sum_{T} (-1)^{|T| - |\lambda|} \prod_{i \ge 1} x_i^{\# \text{$i$'s in $T$}},
$$
where the sum runs over shape $\lambda$ {\it set-valued tableaux} $T$, a generalization of semistandard Young tableaux so that boxes contain sets of integers. 

Being a generalization of the Schur basis, symmetric Grothendieck polynomials share with it many similarities. However $\{ G_{\lambda}\}$ is {\it inhomogeneous} and of {\it unbounded} degree when defined 
for infinitely many variables $(x_1, x_2, \ldots)$, i.e., it is an element of the completion $\hat\Lambda$ of the ring $\Lambda$ of symmetric functions. For example, $G_{(1)} = e_1 - e_2 + e_3 - \cdots$, where $e_k$ is the $k$th elementary symmetric function. It is thus surprising that $\{ G_{\lambda} \}$ is a  {\it basis} of a certain ring: each product $G_{\mu} G_{\nu}$ is a {\it finite} linear combination of $\{G_{\lambda} \}$. When $\mu $ is a single row or column partition, finite expansion was a consequence of Pieri rules proved by Lenart \cite{lenart}. In a general case, this property follows from a Littlewood-Richardson rule given by Buch \cite{buch}. The ring spanned by Grothendieck polynomials is related to the K-theory of Grassmannians \cite{buch}. There is also an important basis $\{g_{\lambda}\}$ of $\Lambda$, dual to $\{G_{\lambda} \}$, that was described via plane partitions and studied by Lam and Pylyavskyy \cite{lp}. %dual to $\{ G_{\lambda}\}$ 

In this paper we study symmetric {\it skew} Grothendieck polynomials via noncommutative {\it Schur operators}. We used these operators in \cite{dy} to prove dualities for certain two-parameter deformations of Grothendieck polynomials. Employing classical Schur operators turns out to be beneficial for obtaining a number of new properties. 

Our main results are the following.

\subsection{Skew Cauchy identity} We prove the following identity %{\it skew Cauchy identity} 
that becomes our central object. 
\begin{theorem}
Let $\mu, \nu$ be any fixed partitions, then\footnote{The shape $\lambda/\!\!/\mu$ is not the usual skew shape $\lambda/\mu$.}
\begin{align*}%\label{skc2}
\sum_{\lambda} G^{}_{\lambda /\!\!/ \mu}(x_1, x_2, \ldots) g^{}_{\lambda/\nu}(y_1, y_2, \ldots) = \prod_{i,j}\frac{1}{1 - x_iy_j} \sum_{\kappa} G^{}_{\nu/\!\!/ \kappa}(x_1, x_2, \ldots) g^{}_{\mu/\kappa}(y_1, y_2, \ldots).
\end{align*}
\end{theorem}
We give a number of applications of {this identity} using it in both operator and generating function forms. Our approach is based on Schur operators as in Fomin~ \cite{fomin}.
For Schur functions such an identity was given by Zelevinsky in the Russian translation of Macdonald's book \cite{macdonald}. Macdonald  \cite[Ch.1 notes]{macdonald} mentioned that this result has apparently been discovered independently by Lascoux, Towber, Stanley, Zelevinsky. It is also known for analogues, e.g., shifted Schur functions \cite{fomin,SS}. Borodin's symmetric functions \cite{bor} generalizing Hall-Littlewood polynomials, also satisfy Cauchy identities which is important in certain stochastic models \cite{bp}; special cases of these symmetric functions have similarity with Grothendieck polynomials \cite{bor}. 

\subsection{Skew Pieri rules} One can immediately obtain some Pieri-type formulas from skew Cauchy identities (e.g., for $\mu = \varnothing$ or $\nu = \varnothing$, which is quite useful). But there is more: it is a general principle that skew Cauchy identities imply {skew} Pieri formulas. For Hall-Littlewood polynomials this was shown by Warnaar \cite{warnaar} using the $q$-binomial theorem for Macdonald polynomials \cite{lw}. In the same way we formulate a general skew Pieri-type formula for dual families of symmetric functions (Lemma \ref{sptf}) and obtain various Pieri formulas for Grothendieck polynomials. We then prove {\it skew Pieri rules} (Theorem~\ref{spieri}) for multiplying the skew Grothendieck polynomials $G_{\mu/\!\!/\nu}$ and $g_{\mu/\nu}$ on elements indexed by single row or column partitions. %For example, in a simplest case we have the following formulas (here $i(\lambda)$ is the number of removable boxes of $\lambda$):
%\begin{align*}
%g_{(1)} g_{\mu/\nu} &= (-i(\mu) + i(\nu)) g_{\mu/\nu} + \sum_{{\lambda = \mu + \square}} g_{\lambda/\nu} - \sum_{{\eta = \nu - \square}} g_{\lambda/\eta}  \\ %\qquad
%G_{(1)} G_{\mu/\!\!/\nu} &= \sum_{{\lambda/\mu \text{ rook strip} \atop \eta \subset \nu}} (-1)^{|\lambda/\mu|} G_{\lambda/\!\!/\eta}.
%\end{align*}
A skew Pieri rule  for Schur functions %(that implies from our general formulas as well) 
was first proved by Assaf and McNamara \cite{AM}. Generalizations were studied in \cite{kl, lls, tw}. 

%We also derive various related generating series. A curious one involves {\it Catalan numbers}:
%$$
%\sum_{\lambda} G_{\lambda/ \delta_{n}} = \text{Cat}_n \prod_{i} \frac{1}{1 - x_i}, \quad \delta_n = (n,n-1,\ldots, 1).
%$$
\subsection{Basis phenomenon} It might seem miraculous that products $G_{\mu} G_{\nu}$ expand {finitely} in the {basis} $\{G_{\lambda}\}$, as known proofs of this property rely on a Littlewood-Richardson (LR) rule  \cite{buch, pp1}, i.e., explicit combinatorial interpretations of multiplicative structure coefficients. 
But is there is a more direct way to see this on the level of symmetric functions? 
We give a quite transparent explanation of this property without appealing to any LR rule. It is what we call a {\it damping} condition of a dual basis is crucial here. Combined with a duality automorphism, the following property of a dual basis prevents infinite expansions: if for a symmetric polynomial %$f_{\lambda/\mu}$ 
we have $f_{\lambda/\mu}(x_1, \ldots, x_n) \ne 0$ then $\lambda_1$ or $\ell(\lambda)$ is bounded from above by a constant depending on $\mu$ and $n$. Symmetric polynomials whose formulas are defined via tableaux with strict row or column conditions have this property. However the dual Grothendieck polynomial $g_{\lambda/\mu}$ does not satisfy it. This issue can be resolved using the polynomials $\omega(g_{\lambda/\mu})$ instead. See Section~\ref{bf} for details and general conditions.    

\subsection{Dual filtered Young graphs} While the homogeneous Schur case corresponds to a graded ring and a combinatorial object behind this is {\it self-dual graded Young's lattice}, Grothendieck polynomials correspond to a filtered ring and {\it dual filtered Young graphs}.  {\it Dual filtered graphs} introduced by Patrias and Pylyavskyy \cite{pp} is a K-theoretic analogue of Stanley's {\it differential posets} \cite{stadiff} and Fomin's {\it dual graded graphs} \cite{fomindual}. Our approach provides new types of {dual filtered Young graphs}. As it was mentioned in \cite{pp}, apparently the most important filtration of dual graded graphs is the so-called {\it M\"obius deformation} as it is related to K-theoretic insertion and LR rules. Even though this deformation (sometimes) produces a dual filtered graph, it is unclear why. We show that the M\"obius deformation of Young's lattice can be obtained from our {\it Cauchy deformation} by a natural transform related to M\"obius inversion. This reveals the presence of a M\"obius deformation on Young's lattice. 
In addition, the constructions of dual filtered graphs give enumerative identities as applications to the normal ordering of differential operators.

\subsection*{Acknowledgements}
I am grateful to Alexei Borodin, Askar Dzhumadil'daev, Thomas Lam, Igor Pak, Leonid Petrov, and Pavlo Pylyavskyy for stimulating and helpful conversations. 

\section{Partitions and Young diagrams}
A {\it partition} is a sequence $\lambda = (\lambda_1 \ge \lambda_2 \ge \ldots)$ of nonnegative integers with only finitely many nonzero terms. The weight of a partition $\lambda$ is the sum $|\lambda| = \lambda_1 + \lambda_2 + \cdots.$ Any partition $\lambda$ is represented as a {\it Young diagram} which contains $\lambda_i$ boxes in its  $i$th row ($i = 1, 2, \ldots$); equivalently it is the set $\{(i, j) : 1 \le i \le \ell, 1 \le j \le \lambda_i \}$, where $\ell = \ell(\lambda)$ is the number of nonzero parts of $\lambda$. The partition $\lambda'$ denotes the {\it conjugate} of $\lambda$ obtained by transposing its diagram. We use English notation for Young diagrams, index columns from left to right and rows from top to bottom.

\ytableausetup{smalltableaux}
The following notation is used throughout the paper. 

\begin{itemize}
\item Let $I(\mu)$ be the set of {\it removable} boxes of $\mu$, it corresponds to {\it inner corners}. Denote 
$i(\mu) := \#I(\mu)$ the number of removable boxes of $\mu$. 
For partitions $\lambda \supset \mu$, define the following extension of a skew shape that will be important: $$\lambda/\!\!/\mu := \lambda/\mu \cup I(\mu).$$ E.g. $(5331)/\!\!/(432)$ consists of the skew shape $(5331)/(432) = \{\scriptsize\ydiagram[*(lightgray)]{1}\}$  and $I(432) = \{\scriptsize\ydiagram[\times]{1} \}$.
%\begin{figure}[t]
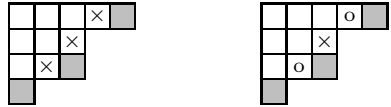
\begin{center}
 \ytableausetup{smalltableaux}
 \begin{ytableau}
~ & ~ & ~ & \times & *(lightgray)\\
~ & ~ & \times \\
~ & \times & *(lightgray) \\
*(lightgray)
\end{ytableau}
%\end{center}
%\begin{center}
\qquad\qquad
\begin{ytableau}
~ & ~ & ~ & \text{o} & *(lightgray)\\
~ & ~ & \times \\
~ & \text{o} & *(lightgray) \\
*(lightgray)
\end{ytableau}
\end{center}
%\caption{Explain}
%\end{figure}
\item Denote by $a({\lambda}/\!\!/{\mu})$ the number of {\it open} boxes (o) of $I(\mu)$ that do not lie in the same column with any box of $\lambda/\mu$. Equivalently, it is just the number of
columns of $\lambda/\!\!/\mu$ that are not columns of $\lambda/\mu$. E.g. $a((5331)/\!\!/(432)) = 2$ (see the diagram above).
\item Denote by $c(\lambda/\mu)$ and $r(\lambda/\mu)$ the number of columns and rows of $\lambda/\mu$. E.g. $c((5331)/(432))=r((5331)/(432))=3$ and $r((5331)/(433)) = 2$. 
%Let $b(\lambda/\mu)$ be the number of connected components of $\lambda/\mu$, e.g., $b((5331)/(43)) = 2$. Let $i(\lambda)$ be the number of {\it inner corners} of $\lambda$, i.e., the number of removable boxes, e.g., $i(5331) = 3$.
\end{itemize}

We also use the following standard terminology: $\lambda/\mu$ is called a {\it horizontal (resp. vertical) strip} if no two boxes of $\lambda/\mu$ lie in the same column (resp. row), equivalently, $|\lambda/\mu| = c(\lambda/\mu)$ (resp. $|\lambda/\mu| = r(\lambda/\mu)$); $\lambda/\mu$ is a {\it rook strip} if no two boxes lie in the same row or column, equivalently, $\lambda/\mu \subset I(\lambda)$.  

Let $\mathbb{Y}$ be the Young lattice, i.e., an infinite graph whose vertices are indexed by partitions and edges are given by $(\lambda, \lambda + \square)$. The {\it M\"obius function} of Young's lattice $\mathbb{Y}$ is given by (e.g., \cite{EC2})
$$
\mu(\lambda, \mu) = 
\begin{cases}
(-1)^{|\lambda/\mu|}, & \text{if $\lambda/\mu$ is a rook strip;}\\
0, & \text{otherwise}.  
\end{cases}
$$
Hence for functions $f, g$ defined on $\mathbb{Y}$, the M\"obius inversion takes the form
\begin{align*}
f(\lambda) = \sum_{\mu \subset \lambda} g(\mu), \qquad g(\lambda) = \sum_{\lambda/\mu \text{ rook strip}} (-1)^{|\lambda/\mu|} f(\mu). 
\end{align*}

\section{Schur operators}
Consider the free $\mathbb{Z}$-module $\mathbb{Z}P = \bigoplus_{\lambda} \mathbb{Z} \cdot \lambda$ with a basis of all partitions. 
\begin{definition}%[Schur operators]
Let ${u} = (u_1, u_2, \ldots)$ and $d = (d_1, d_2, \ldots)$ be sets of linear operators on $\mathbb{Z}P$, called {\it Schur operators}, that act on bases for each $i \ge 1$ as follows: 
\begin{align*}
u_i \cdot \lambda &= \begin{cases}
\lambda + \square\text{ in column $i$}, \text{ if possible,}\\
0, \text{ otherwise}
\end{cases} \quad
%\end{align*}
%\begin{align*}
d_i \cdot \lambda = \begin{cases}
\lambda - \square\text{ in column $i$}, \text{ if possible,}\\
0,  \text{ otherwise.}
\end{cases}
\end{align*}
\end{definition}

The operators $u, d$ build Young diagrams by adding or removing boxes. These operators are noncommutative but they satisfy the following commutation relations that can easily be checked on bases. Let $[a,b] = ab - ba$ denote the commutator. 

\begin{lemma}[\cite{fomin}]\label{lcom1} The following commutation relations hold for the operators $u, d$: 
\begin{itemize}[leftmargin=0.95in]
\item[non-local: ]  
$
[u_j, u_i] = [d_j, d_i] = 0, \ |i- j| \ge 2
$
\item[local Knuth: ]  
$
\label{kcom1} [u_{i+1} u_i, u_i] =  [u_{i+1} u_i, u_{i+1}] =  
[d_{i} d_{i+1}, d_{i}] =  [d_{i} d_{i+1}, d_{i+1}] = 0,  (i \ge 1) 
$
\item[conjugate: ]  
$
[d_j, u_i] = 0\ (i\ne j), \
d_{i+1} u_{i+1} = u_i d_i\ (i \ge 1), \
d_1 u_1 = 1.
$
\end{itemize}
\end{lemma}

Schur operators build Schur polynomials and provide a unified approach for studying their various properties such as Cauchy identities and RSK \cite{fomin}. We generalize this approach for a K-theoretic setting of Grothendieck polynomials.

\subsection{Grothendieck-Schur operators}

Let $\beta$ be a (scalar) parameter and consider the free $\mathbb{Z}[\beta]$-module $\mathbb{Z}[\beta]P = \bigoplus_{\lambda} \mathbb{Z}[\beta] \cdot \lambda$. 

\begin{definition}%[Grothendieck-Schur operators]
Let $\widetilde{u} = (\widetilde{u}_1, \widetilde{u}_2, \ldots)$ and $\widetilde{d} = (\widetilde{d}_1, \widetilde{d}_2, \ldots)$ be linear operators acting on $\mathbb{Z}[\beta]P$ and defined via the Schur operators as follows:
\begin{equation*}
\widetilde{u}_i := u_i - \beta u_i d_i = u_i(1 - \beta d_i), \qquad \widetilde{d}_{i}  := d_i + \beta d_i^2 + \beta^2 d_i^3 + \cdots = (1 - \beta d_i)^{-1} d_i.
\end{equation*}
\end{definition}
For example, \ytableausetup{smalltableaux}
$$\widetilde{u}_2 \cdot {\scriptsize\ydiagram{2, 1}} = {\scriptsize\ydiagram{2,2}} - \beta\ {\scriptsize\ydiagram{2,1}} \qquad \widetilde{u}_2 \cdot {\scriptsize\ydiagram{3,1}} = {\scriptsize\ydiagram{3,2}} \qquad \widetilde{d}_2 \cdot {\scriptsize\ydiagram{3,2,2,2}} = \scriptsize{\ydiagram{3,2,2,1}} +\beta\ {\scriptsize\ydiagram{3,2,1,1}} + \beta^2\ {\scriptsize\ydiagram{3,1,1,1}}$$ 
For a diagram, the operator $\widetilde{u}_i$ adds a box in the $i$th column if possible {\it or} applies the following {\it loop} condition: if the lowest box in the $i$th column is removable it multiplies the result by $-\beta$, since the operator $u_i d_i$ results $1$ (an identity) if the box in the $i$th column is removable, and $0$ otherwise. The operator $\widetilde{d}_i$ removes boxes from the $i$th column graded by $\beta$ in all possible ways. We defined these operators in \cite{dy} when we studied duality properties of stable Grothendieck polynomials.

\begin{lemma}\label{propcom}
The following commutation relations hold for the operators $\widetilde{u}, \widetilde{d}$:
\begin{itemize}[leftmargin=0.85in]
\item[non-local: ]
$
[\widetilde{u}_i, \widetilde{u}_j] = [\widetilde{d}_i, \widetilde{d}_j] = 0,\ |i - j| \ge 2
$
\item[local: ]
$
[\widetilde{u}_{i+1} \widetilde{u}_i, \widetilde{u}_i + \widetilde{u}_{i+1}] =  [\widetilde{d}_{i} \widetilde{d}_{i+1}, \widetilde{d}_i + \widetilde{d}_{i+1}] = 0 \ (i \ge 1)
$
\item[conjugate: ]
$
[\widetilde{u}_i, \widetilde{d}_j] = 0 \  |i - j| \ge 2,\ [\widetilde{u}_{i+1}, \widetilde{d}_i] = 0\ (i \ge 1),\ \widetilde{d}_1 \widetilde{u}_1 = 1.
$ 
\end{itemize}
\end{lemma}
\begin{proof}
In the appendix.
\end{proof}

\begin{remark}
For $\beta = 0$ everything turns into Schur operators. For general $\beta$, the identities are different than in Lemma~\ref{lcom1}. Local Knuth relations do not hold for $\widetilde{u}, \widetilde{d}$, but we have the given local relations instead. Conjugate relations are also different. In general, relations for $\widetilde{u}, \widetilde{d}$ are more complicated than for $u, d$. E.g., the proof of local relations uses the identity $[u_i d_i, u_{i+1}u_i] = 0$ for Schur operators that is not from the list of Lemma \ref{lcom1}. 
\end{remark}

\section{Symmetric skew Grothendieck polynomials}\label{symkg}
Let $x$ be an indeterminate (central variable, commuting with the $u, d$) and define the series
\begin{equation*}
{A}(x) = \cdots (1 + x \widetilde{u}_2)(1 + x \widetilde{u}_1), \qquad {B}(x) = (1 + x \widetilde{d}_1)(1 + x \widetilde{d}_2) \cdots
\end{equation*}
From non-local and local relations given in Lemma \ref{propcom} it is standard (e.g., \cite{fomin,fg}) to deduce that 
\begin{equation*}\label{comab}
[{A}(x), {A}(y)] = 0,\quad [{B}(x), {B}(y)] = 0.
\end{equation*}
Let $\langle \cdot, \cdot \rangle$ be a bilinear pairing on $\mathbb{Z}[\beta]P$ given by $\langle \lambda, \mu \rangle = \delta_{\lambda \mu}.$  % :\mathbb{Z}[\beta]P \times \mathbb{Z}[\beta]P  \to \mathbb{Z}[\beta]
\begin{definition}
Define the {\it skew Grothendieck polynomials} $\{G^{\beta}_{\lambda/\!\!/\mu}\}, \{g^{\beta}_{\lambda/\mu}\}$ as follows
\begin{align*}
G^{\beta}_{\lambda/\!\!/\mu}(x_1, \ldots, x_n) := \langle {A}(x_n) \cdots {A}(x_1) \cdot \mu, \lambda \rangle, \qquad g_{\lambda/\mu}(x_1, \ldots, x_n) := \langle {B}(x_n) \cdots {B}(x_1) \cdot \lambda, \mu \rangle.
%{A}(x_n) \cdots {A}(x_1) \cdot \mu &= \sum_{\lambda} G^{\beta}_{\lambda/\!\!/\mu}(x_1, \ldots, x_n) \cdot \lambda,\quad
%{B}(x_n) \cdots {B}(x_1) \cdot \lambda = \sum_{\mu} g^{\beta}_{\lambda/\mu}(x_1, \ldots, x_n) \cdot \mu. 
\end{align*}
Since the series $[{A}(x_i), A(x_j)] = [B(x_i), B(x_j)] = 0$ commute, the functions $G^{\beta}, g^{\beta}$ (indexed by pairs of partitions) are well-defined polynomials {\it symmetric} in $x_1, \ldots, x_n$. We can then extend these symmetric functions %$G^{\beta}_{\lambda/\!\!/\mu}, g^{\beta}_{\lambda/\mu}$ 
for infinitely many variables $(x_1, x_2, \ldots)$ by letting $n \to \infty$. Notice that $G^{\beta}_{\lambda/\!\!/\mu}= g^{\beta}_{\lambda/\mu} = 0$ if $\mu \not\subseteq \lambda$. 
\end{definition}
\begin{remark}
The reason why we use $G_{\lambda /\!\!/ \mu}$ and {\it not} $G_{\lambda/\mu}$ is in boundary conditions. The function $G_{\lambda/\!\!/\mu}$ depends on the shape $\lambda/\!\!/\mu = \lambda/\mu \cup I(\mu),$ where $I(\mu)$ is the set of removable boxes of $\mu$. For example, we can compute that
$$G^{\beta}_{\lambda/\!\!/\lambda}(x_1, x_2, \ldots) = \prod_{k \ge 1}(1 - \beta x_k)^{i(\lambda)} \ne G_{\varnothing} = 1.$$ 
This is also consistent with the notation in \cite{buch} (though the functions are defined in a different way). 
\end{remark}
\begin{proposition}
The following branching formulas hold:
\begin{align*}
\label{branch}G^{\beta}_{\lambda/\!\!/\mu}(x_1, \ldots, x_n, y_1, \ldots, y_m) &= \sum_{\nu} G^{\beta}_{\lambda/\!\!/\nu}(x_1, \ldots, x_n) G^{\beta}_{\nu/\!\!/\mu}(y_1, \ldots, y_m),\\
g^{\beta}_{\lambda/\mu}(x_1, \ldots, x_n, y_1, \ldots, y_m) &= \sum_{\nu} g^{\beta}_{\lambda/\nu}(x_1, \ldots, x_n) g^{\beta}_{\nu/\mu}(y_1, \ldots, y_m).
\end{align*}
For a single variable $x$ we have
\begin{align*}
G^{\beta}_{\lambda/\!\!/\mu}(x) &= \begin{cases}
(1 - \beta x)^{a({\lambda}/\!\!/{\mu})} x^{|\lambda/\mu|}, & \text{if $\lambda/\mu$ is a horizontal strip;}\\
0, & \text{otherwise} 
\end{cases}\\ 
g^{\beta}_{\lambda/\mu}(x) &= \begin{cases}
\beta^{|\lambda/\mu| - c(\lambda/\mu)}x^{c(\lambda/\mu)}, & \text{if $\mu \subset \lambda$};\\
0, & \text{otherwise.}
\end{cases}
\end{align*}
\end{proposition}
\begin{proof}
The branching formulas are immediate from the definition. For a single variable $x$ we have the following. For any $i$, the operator $1 + x\widetilde{u}_i = 1 + xu_i - x\beta u_i d_i$ applied to $\mu$ can either do nothing, or add a box weighted $x$ in the $i$th column (if possible), or put $-\beta x$ to the removable box of $I(\mu)$ in the $i$th column without growing $\mu$.  The operator $A(x)$ applied to $\mu$ repeats this procedure for $i = 1,2,\ldots$ subsequently. Observe that if some $\widetilde{u}_{{\ell}}$ did not add a box but put $-\beta x$ to the removable box in the column ${\ell}$, then this box will be {\it open}, i.e., there will be no box strictly below it as the operator proceeds with the indices $k > {\ell}$. On the other hand, if $\widetilde{u}_{{\ell}}$ added a new box in the column ${\ell}$, then the last box of $\mu$ in this column will not be open. Furthermore, it is easy to see that the operator $A(x)$ can grow $\mu$ by only horizontal strips. If $\lambda$ is obtained from $\mu$ after applying $A(x)$, then $|\lambda/\mu|$ boxes were added (with the weight $x^{|\lambda/\mu|}$) and {\it some} open boxes from $I(\mu)$ received the weight $-\beta x$. Recall that there are $a(\lambda/\!\!/\mu)$ open boxes. Hence we have
$$
A(x) \cdot \mu %= \prod_{i \ge 1}^{\leftarrow} (xu_i + 1 - x\beta u_i d_i) \cdot \mu = 
= \sum_{\lambda/\mu \text{ hor. strip}} (1 - \beta x)^{a({\lambda}/\!\!/{\mu})} x^{|\lambda/\mu|} \cdot \lambda = \sum_{\lambda/\mu \text{ hor. strip}} G^{\beta}_{\lambda/\!\!/\mu}(x) \cdot \lambda.
$$
The operator $\widetilde{d}_{i_1} \cdots \widetilde{d}_{i_k}$ applied to $\lambda$, removes boxes from the $k$ columns $i_1 < \cdots < i_k$ of $\lambda$ in the order $i_k, \ldots, i_1$ in all possible ways, and thus gives the sum through all $\mu \subset \lambda$  such that $c(\lambda/\mu) = k$ with the corresponding weight $\beta^{|\lambda/\mu| - k}$. Therefore,
$$ 
B(x)\cdot \lambda = \prod_{i \ge 1}^{\rightarrow} (1 + x (d_i + \beta d_i^2 + \cdots)) \cdot \lambda = \sum_{\mu \subset \lambda} \beta^{|\lambda/\mu| - k} x^{c(\lambda/\mu)} \cdot \mu = \sum_{\mu \subset \lambda} g^{\beta}_{\lambda/\mu}(x) \cdot \mu.
$$ 
\end{proof}
In particular, we now obtain the following branching formulas:
\begin{align*}
G^{\beta}_{\lambda/\!\!/\mu}(x_1, \ldots, x_n, x) &= \sum_{\nu/\mu \text{ hor. strip}} G^{\beta}_{\lambda/\!\!/\nu}(x_1, \ldots, x_n) (1 - \beta x)^{a({\nu}/\!\!/{\mu})} x^{|\nu/\mu|} \nonumber
\\ &= \sum_{\lambda/\nu \text{ hor. strip}} (1 - \beta x)^{a({\lambda}/\!\!/{\nu})} x^{|\lambda/\nu|} G^{\beta}_{\nu/\!\!/\mu}(x_1, \ldots, x_n),\\
g^{\beta}_{\lambda/\mu}(x_1, \ldots, x_n, x) &= \sum_{\nu} g^{\beta}_{\lambda/\nu}(x_1, \ldots, x_n) \beta^{|\nu/\mu| - c(\nu/\mu)}x^{c(\nu/\mu)}\nonumber \\
&= \sum_{\nu} \beta^{|\lambda/\nu| - c(\lambda/\mu)}x^{c(\lambda/\nu)} g^{\beta}_{\nu/\mu}(x_1, \ldots, x_n). 
\end{align*}

\begin{figure}
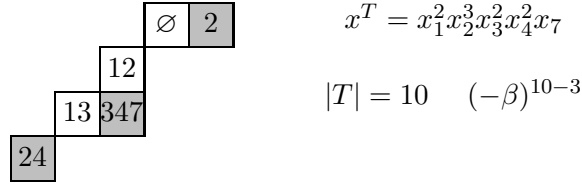

\begin{center}
{%\scriptsize
 \ytableausetup{boxsize=normal}
 \begin{ytableau}
\none & \none & \none & \varnothing & *(lightgray) 2\\
\none & \none & 12 \\
\none & 13 & *(lightgray) 347 \\
*(lightgray) 24
\end{ytableau}
}
\qquad 
\begin{tabular}{c}
\\
$x^T = x_1^2 x_2^3 x_3^2 x_4^2 x_7$ \\ 
\\
$|T| = 10$ \quad $(-\beta)^{10 - 3}$
\end{tabular}
\end{center}
\caption{A set-valued tableau of shape $(5331)/\!\!/(432) = (5331)/(432) \cup I(432)$. The gray boxes form the skew shape $(5331)/(432)$; the white cells are the removable boxes $I(432)$ of $(432)$. It is allowed to put $\varnothing$ in $I(\mu)$ but boxes of $\lambda/\mu$ must be nonempty.}
\label{svtfig}
\end{figure}

\begin{definition}
A {\it set-valued tableaux} (SVT) of shape $\lambda/\!\!/\mu$ is a filling of boxes of $\lambda/\!\!/\mu = \lambda/\mu \cup I(\mu)$ (skew shape $\lambda/\mu$ and removable boxes of $\mu$) by sets of positive integers such that if one replaces each set by any of its elements the resulting tableau is {\it semistandard} (i.e., has weakly increasing rows from left to right and strictly increasing columns from top to bottom). When filling $\lambda/\mu$, sets in boxes should be {\it nonempty}, however when filling the boxes of $I(\mu)$ it is {allowed} to have {\it empty} sets. 
For a set-valued tableau $T$, the corresponding monomial is defined as $x^T = \prod_{i \ge 1} x_i^{a_i}$, where $a_i$ is the number of $i$'s in $T$ and let $|T| = \sum_{i} a_i$. See Fig.~\ref{svtfig}.
\end{definition}
\begin{definition}
A {\it reverse plane partition} (RPP) of shape $\lambda/\mu$ is a filling of a Young diagram of $\lambda/\mu$ by positive integers weakly increasing in rows from left to right and columns from top to bottom. For a reverse plane partition $T$, the corresponding monomial is defined as $x^T = \prod_{i \ge 1} x_i^{c_i}$ where $c_i$ is the number of columns of $T$ containing $i$ and let $|T| = \sum_{i} c_i$.
\end{definition}
\begin{theorem}The following combinatorial formulas hold: 
\begin{align*}
G^{\beta}_{\lambda/\!\!/\mu} &= \sum_{T \in SVT(\lambda/\!\!/\mu)} (-\beta)^{|T| - |\lambda/\mu|} x^T,\qquad g^{\beta}_{\lambda/\mu} = \sum_{T \in RPP(\lambda/\mu)} \beta^{|\lambda/\mu| - |T|} x^T.
\end{align*}
\end{theorem}
\begin{proof}
We construct any set-valued tableau recursively by applying the operators $A(x_1), A(x_2), \ldots$ and recording the entries $1,2, \ldots$ in diagrams. Applying the factor $(1 + x_k u_i - x_k \beta u_i d_i)$ from the operator $A(x_k)$ to the current tableau, we either do nothing or add a new box in the $i$th column and record the entry $k$ in this box (notice that this box is removable and will be treated so further), or if the box in the $i$th column is removable we add $k$ to the existing entries of this box. It is clear that semistandard inequalities are preserved during these operations. The formula for $g^{\beta}_{\lambda/\mu}$ can be explained similarly or it is easily seen by combining the branching and single variable formulas. 
\end{proof}
\begin{remark}
Using the the operators $\widetilde{u}, \widetilde{d}$ we obtained combinatorial formulas due to Buch \cite{buch} and Lam-Pylyavskyy \cite{lp}.
For $\mu = \varnothing$ and $\beta = 1$ they coincide with dual families of stable Grothendieck polynomials $\{G_{\lambda}\}, \{g_{\lambda}\}$. They are Hopf-dual or dual via the Hall inner product $\langle G_{\lambda}, g_{\mu}\rangle = \delta_{\lambda \mu}$ for which Schur functions form an orthonormal basis.   
\end{remark}
\begin{remark}
Let $\beta = 1$ and $G_{\lambda/\mu}$ be the symmetric function defined via SVT formula of skew shape $\lambda/\mu$ (without the removable boxes $I(\mu)$ of $\mu$). Then the two functions are related \cite{buch} via the M\"obius inversion  
\begin{align*}
G_{\lambda/\mu} = \sum_{\nu \subset \mu} G_{\lambda/\!\!/\nu}, \qquad G_{\lambda/\!\!/\nu} = \sum_{\nu/\mu \text{ rook strip}} (-1)^{|\nu/\mu|} G_{\lambda/\mu}. 
\end{align*}
\end{remark}

\section{Skew Cauchy identity}
We use the notation $\mathbf{x} = (x_1, x_2, \ldots)$ and $\mathbf{y} = (y_1, y_2, \ldots)$.
\begin{theorem}%[Skew Cauchy identity] 
\label{skc}
Let $\mu, \nu$ be any fixed partitions, then
\begin{align}%\label{skc2}
\sum_{\lambda} G^{\beta}_{\lambda /\!\!/ \mu}(\mathbf{x}) g^{\beta}_{\lambda/\nu}(\mathbf{y}) = \prod_{i,j}\frac{1}{1 - x_iy_j} \sum_{\kappa} G^{\beta}_{\nu/\!\!/ \kappa}(\mathbf{x}) g^{\beta}_{\mu/\kappa}(\mathbf{y}).
\end{align}
%\begin{align}%\label{skc2}
%\sum_{\lambda} G^{\beta}_{\lambda /\!\!/ \mu}(x_1, x_2, \ldots) g^{\beta}_{\lambda/\nu}(y_1, y_2, \ldots) = \prod_{i,j}\frac{1}{1 - x_iy_j} \sum_{\kappa} G^{\beta}_{\nu/\!\!/ \kappa}(x_1, x_2, \ldots) g^{\beta}_{\mu/\kappa}(y_1, y_2, \ldots).
%\end{align}
\end{theorem}
The identity is equivalent to the following commutation relation %(Yang-Baxter-type) 
for the series ${A}, {B}$ defined in Sec.~\ref{symkg}. %For $\beta = 0$ this identity gives the skew Cauchy identity for Schur polynomials.

\begin{theorem}\label{cauchy1} The following commutation relation holds  
\begin{equation*}\label{gcauchy}
{B}(y) {A}(x) =  \frac{1}{1 - xy} {A}(x) {B}(y).
\end{equation*}
\end{theorem}

The proof %of Theorem \ref{cauchy1} 
uses the following Yang-Baxter-type local identity for the operators $\widetilde{u}, \widetilde{d}$.
\begin{lemma}\label{locud} For all $i \ge 1$ we have
$$
(1 - xy \widetilde{u}_i \widetilde{d}_i)^{-1} (1 + x\widetilde{u}_i)(1+y \widetilde{d}_{i+1}) = 
(1 - xy \widetilde{d}_{i+1} \widetilde{u}_{i+1})^{-1} (1 + y\widetilde{d}_{i+1})(1+x \widetilde{u}_{i}).  \eqno{(*)}
$$
\end{lemma}
\begin{proof}
See the Appendix.
\end{proof}

\begin{proof}[Proof of Theorem \ref{cauchy1}]
Using Lemmas \ref{propcom}, \ref{locud} and the identity
$$
(1 + a)(1 - ba)^{-1}(1 + b) = (1 + b)(1 - ab)^{-1}(1+a) \eqno{(**)}
$$
holding for any non-commuting $a$ and $b$, we have
{%\small
\begin{align*}
A(x) (1 - xy)^{-1} & B(y) = \cdots (1 + x \widetilde{u}_2)\underbrace{(1 + x \widetilde{u}_1) (1 - xy \widetilde{d}_1 \widetilde{u}_1)^{-1} (1 + y \widetilde{d}_1)}_{(**)} (1 + y \widetilde{d}_2) \cdots \\
&=  \cdots (1 + x \widetilde{u}_2) \underbrace{(1 + y \widetilde{d}_1)}_{\leftarrow} \ \underbrace{(1 - xy \widetilde{u}_1 \widetilde{d}_1)^{-1} (1 + x \widetilde{u}_1) (1 + y \widetilde{d}_2)}_{(*)} \cdots \\
&= (1 + y \widetilde{d}_1) \cdots (1 + x \widetilde{u}_3)(1 + x \widetilde{u}_2) \ (1 - xy \widetilde{d}_2 \widetilde{u}_2)^{-1} (1 + y \widetilde{d}_2) \underbrace{(1 + x \widetilde{u}_1)}_{\rightarrow}(1 + y\widetilde{d}_3) \cdots \\
&= (1 + y \widetilde{d}_1) \cdots (1 + x \widetilde{u}_3)\underbrace{(1 + x \widetilde{u}_2)  (1 - xy \widetilde{d}_2 \widetilde{u}_2)^{-1} (1 + y \widetilde{d}_2)}_{(**)} (1 + y \widetilde{d}_3) \cdots (1 + x \widetilde{u}_1) \\
& \cdots \\
&= B(y) A(x).
\end{align*}
}
\end{proof}

\begin{proof}[Proof of Theorem \ref{skc}]
Using Theorem \ref{cauchy1} and the commutativity $[A(x_i), A(x_j)] = [B(y_i), B(y_j)] = 0$ we have
\begin{align*}
\sum_{\lambda} G^{\beta}_{\lambda /\!\!/ \mu}(\mathbf{x}) g^{\beta}_{\lambda/\nu}(\mathbf{y})
&= \langle  \prod_{j} B(y_j) \prod_{i} A(x_i) \cdot \mu, \nu \rangle \\
&= \prod_{i,j}\frac{1}{1 - x_i y_j} \langle  \prod_{i} A(x_i) \prod_{j} B(y_j) \cdot \mu, \nu \rangle \\
&= \prod_{i,j}\frac{1}{1 - x_i y_j} \sum_{\kappa} G^{\beta}_{\nu/\!\!/ \kappa}(\mathbf{x}) g^{\beta}_{\mu/\kappa}(\mathbf{y}).
\end{align*}
\end{proof}

\begin{remark}
In fact, Theorem \ref{cauchy1} is equivalent to the {\it single} variable Cauchy identity
\begin{equation*}%\label{grone}
\sum_{\lambda} G^{\beta}_{\lambda/\!\!/\mu}(x) g^{\beta}_{\lambda/\nu}(y) = \frac{1}{1 - xy} \sum_{\kappa} G^{\beta}_{\nu/\!\!/\kappa}(x) g^{\beta}_{\mu/\kappa}(y).
\end{equation*}
So to prove Theorem \ref{skc} it is enough to prove it just for {single} variables $x,y$.
\end{remark}

\subsection{Corollaries}
First note that by setting $\beta = 0$, the results generalize corresponding properties of skew Schur polynomials, since $G^{0}_{\lambda /\!\!/ \mu} = g^{0}_{\lambda/\mu} = s_{\lambda/\mu}.$

\begin{corollary}[Cauchy identity]
For $\mu = \nu = \varnothing$ we obtain the usual Cauchy identity
\begin{align*}
\sum_{\lambda} G^{\beta}_{\lambda}(\mathbf{x}) g^{\beta}_{\lambda}(\mathbf{y}) = \prod_{i,j} \frac{1}{1 - x_i y_j}.
\end{align*}
This identity is equivalent to the duality of families $\{G_{\lambda} \}, \{g_{\mu} \}$ via the standard Hall inner product for which Schur functions are orthonormal, i.e., $\langle G_{\lambda}, g_{\mu}\rangle = \delta_{\lambda \mu}$. 
\end{corollary}

\begin{corollary}[Pieri-type formulas]
For $\mu = \varnothing$ or $\nu = \varnothing$ we have
\begin{align}%\label{pieritype}
        \sum_{\lambda} G^{\beta}_{\lambda}(\mathbf{x}) g^{\beta}_{\lambda/\nu}(\mathbf{y}) &= \prod_{i,j}\frac{1}{1 - x_iy_j} G^{\beta}_{\nu}(\mathbf{x}),\label{pieritype1}\\
        \sum_{\lambda} G^{\beta}_{\lambda /\!\!/ \mu}(\mathbf{x}) g^{\beta}_{\lambda}(\mathbf{y}) &= \prod_{i,j}\frac{1}{1 - x_iy_j} g^{\beta}_{\mu}(\mathbf{y}).\label{pieritype2}
\end{align}
\end{corollary}
In particular, the last identity gives the following formulas for $\mathbf{y} = (1,0,0,\ldots)$.
\begin{align*}
\sum_{\lambda}  \beta^{|\lambda| - c(\lambda)} G^{\beta}_{\lambda/\!\!/\mu} &= \beta^{|\mu| - c(\mu)}\prod_{i} \frac{1}{1 - x_i}, \qquad  \sum_{\lambda}  G_{\lambda/\!\!/\mu} = \prod_{i} \frac{1}{1 - x_i}.
\end{align*}
Recall in contrast a similar identity for Schur functions (e.g., \cite{macdonald}):
\begin{align*}
\sum_{\lambda} s_{\lambda/\mu} = \prod_{i}\frac{1}{1 - x_i} \prod_{i < j} \frac{1}{1 - x_i x_j} \sum_{\kappa} s_{\mu/\kappa}.
\end{align*}
%Let $\mathbf{y} = (y_1, y_2, \ldots)$. 
The Pieri-type formulas specialize to the following identities for $\mathbf{y} = (q,0,0, \ldots)$
\begin{align*}
    \sum_{\lambda}  q^{c(\lambda/\nu)} \beta^{|\lambda/\nu| - c(\lambda/\nu)} G^{\beta}_{\lambda} &= \prod_{i} \frac{1}{1 - q x_i} G^{\beta}_{\nu},\\
    \sum_{\lambda/\mu \text{ hor. strip}} (1 - \beta q)^{a(\lambda/\!\!/\mu)} q^{|\lambda/\mu|} g^{\beta}_{\lambda} &= \prod_{i} \frac{1}{1 - q x_i} g^{\beta}_{\mu}.
\end{align*}
Consider another specialization:
\begin{align*}
\sum_{\lambda} G^{\beta}_{\lambda}(\mathbf{x}) g^{\beta, q}_{\lambda}(\mathbf{y}) &= \prod_{i,j} \frac{1}{1 - x_i y_j} \prod_{i} \frac{1}{1 - q x_i},\\
g^{\beta, q}_{\lambda} &:= g^{\beta}_{\lambda}(q,\mathbf{y}) = \sum_{\nu \subset \lambda} q^{c(\lambda/\nu)} \beta^{|\lambda/\nu| - c(\lambda/\nu)} g^{\beta}_{\lambda/\nu}(\mathbf{y})
\end{align*}
Let $d(\lambda) := \#\{\mu : \mu \subset \lambda \}$ be the number of subdiagrams of $\lambda$. Set $\mathbf{y} = (1,0,\ldots)$ and $\beta = q = 1$: %the last identity gives the following
$$
\sum_{\lambda} d(\lambda) G_{\lambda} = \prod_{i} \frac{1}{(1 - x_i)^2}, %\qquad d(\lambda) := \#\{\mu : \mu \subset \lambda \}.
$$
For a `pure' skew shape $\lambda/\mu$ we obtain the following generating series.
\begin{proposition}
We have
\begin{align*}
\sum_{\lambda} G_{\lambda/\mu}(\mathbf{x}) g_{\lambda}(\mathbf{y}) &= \prod_{i,j} \frac{1}{1 - x_i y_j} g_{\mu}(1, \mathbf{y}).
\end{align*}
\end{proposition}
\begin{proof}
Recall that 
$G_{\lambda/\mu} = \sum_{\nu \subset \mu} G_{\lambda/\!\!/\nu}.$
From \eqref{pieritype2} we have
\begin{align*}
\prod_{i,j} \frac{1}{1 - x_i y_j} g_{\mu}(1,\mathbf{y}) &= \prod_{i,j} \frac{1}{1 - x_i y_j} \sum_{\nu \subset \mu} g_{\nu}(\mathbf{y}) \\
&= \sum_{\lambda} \sum_{\nu \subset \mu} G_{\lambda/\!\!/\nu}(\mathbf{x}) g_{\lambda}(\mathbf{y}) \\
&= \sum_{\lambda} G_{\lambda/\mu}(\mathbf{x}) g_{\lambda}(\mathbf{y}).
\end{align*}
\end{proof}
This gives the next formulas including a curious identity involving the {\it Catalan numbers} $\text{Cat}_n = \binom{2n}{n} / (n+1)$.
\begin{corollary}
\begin{align*} 
\sum_{\lambda}  G_{\lambda/\mu} &= d(\mu) \prod_{i} \frac{1}{1 - x_i},\qquad %\qquad \sum_{\lambda} d(\lambda) G_{\lambda} = \prod_{i} \frac{1}{(1 - x_i)^2}
\sum_{\lambda}  G_{\lambda/\delta_n} = \text{Cat}_n \prod_{i} \frac{1}{1 - x_i}, \quad \delta_{n} = (n, n-1, \ldots, 1).
\end{align*}
\end{corollary}

\section{Dual functions via the standard involution}\label{jjj}
We also need to describe the following dual families of symmetric functions. Let
\begin{equation*}
\overline{A}(x) = (1 - x \widetilde{u}_1)^{-1} (1 - x \widetilde{u}_2)^{-1} \cdots, \qquad \overline{B}(x) =  \cdots(1 - x \widetilde{d}_2)^{-1}(1 - x \widetilde{d}_1)^{-1}
\end{equation*}
so that $A(x) \overline{A}(-x) = B(x)\overline{B}(-x) = 1$.
Again, using non-local and local commutativity given in Lemma~\ref{propcom} we deduce that  
\begin{equation*}%\label{comab}
[\overline{A}(x), \overline{A}(y)] = 0,\quad [\overline{B}(x), \overline{B}(y)] = 0.
\end{equation*}
Hence we can define the symmetric functions $\{J^{\beta}_{\lambda/\!\!/\mu}\}, \{j^{\beta}_{\lambda/\mu}\}$ via the series  
\begin{align*}\label{jjdef}
J^{\beta}_{\lambda/\!\!/\mu}(x_1, \ldots, x_n) := \langle \overline{A}(x_n) \cdots \overline{A}(x_1) \cdot \mu, \lambda \rangle, \qquad j^{\beta}_{\lambda/\mu}(x_1, \ldots, x_n) := \langle \overline{B}(x_n) \cdots \overline{B}(x_1) \cdot \lambda, \mu \rangle.
%\overline{A}(x_n) \cdots \overline{A}(x_1) \cdot \mu &= \sum_{\lambda} J^{\beta}_{\lambda/\!\!/\mu}(x_1, \ldots, x_n) \cdot \lambda,\quad
%\overline{B}(x_n) \cdots \overline{B}(x_1) \cdot \lambda = \sum_{\mu} j^{\beta}_{\lambda/\mu}(x_1, \ldots, x_n) \cdot \mu. 
\end{align*}
As for $G, g$, it is not hard to obtain the following formulas for a single variable $x$ 
\begin{align}
J^{\beta}_{\lambda/\!\!/\mu}(x) &= 
	\begin{cases}
		\displaystyle\frac{1}{(1+ \beta x)^{a(\lambda'/\!\!/\mu')}} \left(\frac{x}{1 + \beta x}\right)^{|\lambda/\mu|}, & \text{ if $\lambda/\mu$ is a vert. strip},\\
		0, & \text{ otherwise,}
	\end{cases}\\
j^{\beta}_{\lambda/\!\!/\mu}(x) &= 
	\begin{cases}
		x^{c(\lambda/\mu)}(\beta + x)^{|\lambda/\mu| - c(\lambda/\mu)}, & \text{ if $\lambda/\mu$ is a vert. strip},\\
		0, & \text{ otherwise.}
	\end{cases}
\end{align}
Define {\it multiset-valued tableaux} (MSVT) of shape $\lambda/\!\!/\mu$ to be the filling as for set-valued tableaux %(see Sec. \ref{?}) 
but allowing multisets (i.e., sets with repeated elements) in boxes. The corresponding monomial is defined in the same way: $x^{T} = \sum_{i} x_i^{a_i},$ where $a_i$ is the number of $i$'s in $T$ and $|T| = \sum_i a_i$. For a  semistandard Young tableaux (SSYT) $T$, let $r_i$ be the number of rows containing $i$ and $a_i$ be the total number of $i$'s. 
\begin{proposition}
The following combinatorial formulas hold
\begin{align*}
J^{\beta}_{\lambda/\!\!/\mu} = \sum_{T \in MSVT(\lambda'/\!\!/\mu')} (-\beta)^{|T| - |\lambda/\mu|} x^T, \qquad j^{\beta}_{\lambda/\mu} = \sum_{T \in SSYT(\lambda'/\mu')} \prod_{i \in T} x_i^{r_i}(x_i + \beta)^{a_i - r_i}.
\end{align*}
\end{proposition}

Let $\omega$ be the standard ring automorphism satisfying $\omega : s_{\lambda} \mapsto s_{\lambda'}, $ where $\lambda'$ is the conjugate partition. 
To establish the next result we use the method of Fomin and Greene \cite{fg} of {\it noncommutative Schur functions}. Let $\mathbf{\widetilde{u}} = (\widetilde{u}_1, \widetilde{u}_2, \ldots)$ and define the noncommutative Schur functions as follows
\begin{align*}
s_{\lambda}(\mathbf{\widetilde{u}}) &:= \det[e_{\lambda' - i + j}(\widetilde{\mathbf{u}})] = \sum_{\sigma \in S_{\lambda_1}} \sgn(\sigma) e_{\lambda'_1 + \sigma(1) - 1}(\widetilde{\mathbf{u}}) \cdots e_{\lambda'_1 + \sigma(\lambda_1) - \lambda_1}(\widetilde{\mathbf{u}}),\\
e_k(\widetilde{\mathbf{u}}) &:= \sum_{i_1 > \cdots > i_k \ge 1} \widetilde{u}_{i_1} \cdots \widetilde{u}_{i_k}.
\end{align*} 
Since the operators $\widetilde{u}$ satisfy non-local and local commutativity relations (Lemma \ref{propcom})
$$
[\widetilde{u}_i, \widetilde{u}_j] = 0,\ |i - j| \ge 2, \qquad [\widetilde{u}_{i+1}\widetilde{u}_i, \widetilde{u}_{i} + \widetilde{u}_{i+1}] = 0,
$$
we conclude that these functions commute:
$$
[e_i(\mathbf{\widetilde{u}}), e_j(\mathbf{\widetilde{u}})] = [s_{\lambda}(\mathbf{\widetilde{u}}), s_{\mu}(\mathbf{\widetilde{u}})] = 0,\quad \forall i,j, \lambda, \mu
$$
and the following noncommutative versions of Cauchy identities hold:
\begin{align*}
\cdots A(x_2) A(x_1) &=\prod_{i} \prod_{j}^{\leftarrow} (1 + x_i \widetilde{u}_j) =\sum_{\lambda} s_{\lambda'}(\mathbf{x}) s_{\lambda}(\mathbf{\widetilde{u}}),\\
\cdots \overline{A}(x_2) \overline{A}(x_1) &= \prod_{i} \prod_{j}^{\rightarrow} (1 - x_i \widetilde{u}_j)^{-1} =\sum_{\lambda} s_{\lambda}(\mathbf{x}) s_{\lambda}(\mathbf{\widetilde{u}})
\end{align*}
The same holds for the operators $\widetilde{d}$ and the corresponding series $B(x), \overline{B}(x)$. 
\begin{theorem}\label{omega} We have %The following properties hold: 
$\omega(G^{\beta}_{\lambda/\!\!/\mu}) = J^{\beta}_{\lambda/\!\!/\mu}$ and $\omega(g^{\beta}_{\lambda/\mu}) = j^{\beta}_{\lambda/\mu}$.
\end{theorem}
\begin{proof}
Using noncommutative Schur functions and Cauchy identities we have
\begin{align*}
\omega(G^{\beta}_{\lambda/\!\!/\mu}) 
	&= \omega \langle \prod_i A(x_i) \cdot \mu, \lambda \rangle \\
	&= \omega \langle \sum_{\nu} s_{\nu'}(\mathbf{x}) s_{\nu}(\mathbf{\widetilde{u}}) \cdot \mu, \lambda\rangle\\
	&= \sum_{\nu} \omega(s_{\nu'}(\mathbf{x})) \langle s_{\nu}(\mathbf{\widetilde{u}}) \cdot \mu, \lambda\rangle\\
	&= \sum_{\nu} s_{\nu}(\mathbf{x}) \langle s_{\nu}(\mathbf{\widetilde{u}}) \cdot \mu, \lambda\rangle	\\
	&= \langle \sum_{\nu} s_{\nu}(\mathbf{x}) s_{\nu}(\mathbf{\widetilde{u}}) \cdot \mu, \lambda\rangle	\\
	&= \langle \prod_{i} \overline{A}(x_i) \cdot \mu, \lambda \rangle \\
	&= J^{\beta}_{\lambda/\!\!/\mu}.
\end{align*}
The property $\omega(g^{\beta}_{\lambda/\mu}) = j^{\beta}_{\lambda/\mu}$ follows in the same way.
\end{proof}

\begin{corollary}[Dual skew Cauchy identities] We have
\begin{align*}
\sum_{\lambda} J^{\beta}_{\lambda /\!\!/ \mu}(\mathbf{x}) g^{\beta}_{\lambda/\nu}(\mathbf{y}) &= \prod_{i,j}{(1 + x_iy_j)} \sum_{\kappa} J^{\beta}_{\nu/\!\!/ \kappa}(\mathbf{x}) g^{\beta}_{\mu/\kappa}(\mathbf{y}),\\
\sum_{\lambda} G^{\beta}_{\lambda /\!\!/ \mu}(\mathbf{x}) j^{\beta}_{\lambda/\nu}(\mathbf{y}) &= \prod_{i,j}{(1 + x_iy_j)} \sum_{\kappa} G^{\beta}_{\nu/\!\!/ \kappa}(\mathbf{x}) j^{\beta}_{\mu/\kappa}(\mathbf{y}),\\
\sum_{\lambda} J^{\beta}_{\lambda /\!\!/ \mu}(\mathbf{x}) j^{\beta}_{\lambda/\nu}(\mathbf{y}) &= \prod_{i,j}\frac{1}{(1 - x_iy_j)} \sum_{\kappa} J^{\beta}_{\nu/\!\!/ \kappa}(\mathbf{x}) j^{\beta}_{\mu/\kappa}(\mathbf{y}).
\end{align*}
\end{corollary}

\begin{remark}
The functions $J, j$ were first introduced and studied in \cite{lp} using noncommutative operators that are different to ours.
\end{remark}

\begin{remark}
It can be proved that {\it canonical} deformations of symmetric Grothendieck polynomials studied in \cite{dy}, the functions  $G^{(\alpha, \beta)}_{\lambda},$ $g^{(\alpha, \beta)}_{\lambda},$ satisfying $\omega (G^{(\alpha, \beta)}_{\lambda}) = G^{(\beta, \alpha)}_{\lambda'},$ $\omega (g^{(\alpha, \beta)}_{\lambda}) = g^{(\beta, \alpha)}_{\lambda'}$ and $G^{(0, \beta)}_{\lambda} = G^{\beta}_{\lambda}, g^{(0, \beta)}_{\lambda} = g^{\beta}_{\lambda}$, $G^{(\beta, 0)}_{\lambda} = J^{\beta}_{\lambda'}, g^{(\beta, 0)}_{\lambda} = j^{\beta}_{\lambda'}$, also satisfy skew Cauchy identities.
\end{remark}

\section{Skew Pieri formulas}
%Skew Cauchy identity can be used to obtain skew Pieri formulas.  
\begin{theorem}[Skew Pieri rules]\label{spieri} The following formulas hold
\begin{align*}
G^{\beta}_{(1^k)} G^{\beta}_{\mu/\!\!/\nu} &= \sum_{{\lambda/\mu \text{ vert strip} \atop \eta \subset \nu}} W^{\lambda, \mu}_{\nu, \eta}\ G^{\beta}_{\lambda/\!\!/\eta} \quad W^{\lambda, \mu}_{\nu, \eta} = (-1)^{|\lambda/\mu| - k} \beta^{|\lambda/\mu| + |\nu/\eta| - k} \text{\scriptsize $\binom{c(\lambda/\mu) + c(\nu/\eta) - 1}{|\lambda/\mu| + c(\nu/\eta) - k}$} \\ %since vert strip b(\lambda/\mu) = c(\lambda/\mu)
g^{\beta}_{(k)} g^{\beta}_{\mu/\nu} &= \sum_{ {\lambda/\mu \text{ hor strip} \atop \nu/\eta \text{ vert strip}} } w^{\lambda, \mu}_{\nu, \eta}\ g^{\beta}_{\lambda/\eta} \qquad w^{\lambda, \mu}_{\nu, \eta} =  (-1)^{k - |\lambda/\mu|} \beta^{k - |\lambda/\mu| - |\nu/\eta|}  \text{\scriptsize $\binom{a(\lambda/\!\!/\mu) - a(\nu'/\!\!/\eta') - |\nu/\eta|}{k - |\lambda/\mu| - |\nu/\eta|}$} 
\end{align*}
\end{theorem}
\begin{remark}
Note that these expansions are finite. For $\nu = \varnothing$ and $G_{\mu}$ we recover Pieri formulas proved in \cite{lenart}. 
One could interpret the coefficients $w, W$ in various ways: as a number of certain tableaux or as a number of walks in dual graphs defined later. 
\end{remark}

Consider the automorphisms $\widehat\tau : G^{\beta}_{\lambda} \mapsto G^{\beta}_{\lambda'}$ and $\tau : g^{\beta}_{\lambda} \mapsto g^{\beta}_{\lambda'}$. (Such automorphisms exist \cite{buch, dy}.) 
\begin{proposition}
We have $\widehat\tau(G^{\beta}_{\lambda/\!\!/\mu}) = G^{\beta}_{\lambda'/\!\!/\mu'}$ and  $\tau(g^{\beta}_{\lambda/\mu}) = g^{\beta}_{\lambda'/\mu'}$. %and hence the formulas give dual skew Pieri rules as well.
\end{proposition}
\begin{proof}
From the Cauchy identity we have 
$$
\tau_{\mathbf{y}}\widehat{\tau}_{\mathbf{x}} \prod_{i,j}\frac{1}{1 - x_i y_j} = \sum_{\lambda} G^{\beta}_{\lambda'}(\mathbf{x}) g^{\beta}_{\lambda'}(\mathbf{y}) = \prod_{i,j} \frac{1}{1 - x_i y_j}.
$$
Applying ${\tau}_{\mathbf{y}} \widehat{\tau}_{\mathbf{x}}$ on the Pieri-type formula \eqref{pieritype1} and then by using it again we obtain
\begin{align*}
\sum_{\lambda} G^{\beta}_{\lambda'}(\mathbf{x}) \tau_{\mathbf{y}}( g^{\beta}_{\lambda/\nu}(\mathbf{y}) ) &= {\tau}_{\mathbf{y}}\widehat{\tau}_{\mathbf{x}}\left(\prod_{i,j} \frac{1}{1 - x_i y_j} \right) G^{\beta}_{\nu'}(\mathbf{x}) \\
&= \prod_{i,j} \frac{1}{1 - x_i y_j} G^{\beta}_{\nu'}(\mathbf{x}) \\
&= \sum_{\lambda} G^{\beta}_{\lambda}(\mathbf{x}) g^{\beta}_{\lambda/\nu'}(\mathbf{y})
\end{align*}
from which we conclude that $\tau_{\mathbf{y}}( g^{\beta}_{\lambda/\nu}(\mathbf{y}) ) = g^{\beta}_{\lambda'/\nu'}(\mathbf{y})$. The proof of the first formula is the same.
\end{proof}

\begin{corollary}[Dual skew Pieri rules] Applying $\widehat\tau, \tau$ we obtain
\begin{align*}
G^{\beta}_{(k)} G^{\beta}_{\mu/\!\!/\nu} &= \sum_{{\lambda/\mu \text{ hor strip} \atop \eta \subset \nu}} V^{\lambda, \mu}_{\nu, \eta}\ G^{\beta}_{\lambda/\!\!/\eta} \quad V^{\lambda, \mu}_{\nu, \eta} = (-1)^{|\lambda/\mu| - k} \beta^{|\lambda/\mu| + |\nu/\eta| - k} \text{\scriptsize $\binom{r(\lambda/\mu) + r(\nu/\eta) - 1}{|\lambda/\mu| + r(\nu/\eta) - k}$} \\ %since vert strip b(\lambda/\mu) = c(\lambda/\mu)
g^{\beta}_{(k)} g^{\beta}_{\mu/\nu} &= \sum_{ {\lambda/\mu \text{ vert strip} \atop \nu/\eta \text{ hor strip}} } v^{\lambda, \mu}_{\nu, \eta}\ g^{\beta}_{\lambda/\eta} \qquad v^{\lambda, \mu}_{\nu, \eta} = (-1)^{k - |\lambda/\mu|}\beta^{k - |\lambda/\mu| - |\nu/\eta|} \text{\scriptsize $\binom{a(\lambda'/\!\!/\mu') - a(\nu/\!\!/\eta) - |\nu/\eta|}{k - |\lambda/\mu| - |\nu/\eta|}$} 
\end{align*}
\end{corollary}

\begin{corollary}[Simple skew Pieri rules]
Let $\beta = 1$.
We have
\begin{align*}
g_{(1)} g_{\mu/\nu} &= (-i(\mu) + i(\nu)) g_{\mu/\nu} + \sum_{{\lambda = \mu + \square}} g_{\lambda/\nu} - \sum_{{\eta = \nu - \square}} g_{\lambda/\eta}  \\ %\qquad
G_{(1)} G_{\mu/\!\!/\nu} &= \sum_{{\lambda/\mu \text{ rook strip} \atop \eta \subset \nu}} (-1)^{|\lambda/\mu|} G_{\lambda/\!\!/\eta}.
\end{align*}
\end{corollary}
\begin{corollary}
$\beta = 0$ gives the skew Pieri rule for Schur functions $s_{\lambda/\mu}$ \cite{AM}
 $$
s_{(k)} s_{\mu/\nu} = \sum_{ {\lambda/\mu \text{ hor strip} \atop \nu/\eta \text{ vert strip}} } (-1)^{|\nu/\eta|} s_{\lambda/\eta}, \qquad |\lambda/\mu| +|\nu/\eta| = k.
$$
\end{corollary}

\subsection{Dual skew families}
\begin{definition}
    The families $\{F_{\lambda/\mu} \},$ $\{f_{\lambda/\mu} \}$ of symmetric functions (generally lying in the completion $\hat\Lambda$ of the ring of symmetric functions) indexed by pairs of partitions (with boundary conditions $f_{\varnothing/\mu} = \delta_{\varnothing, \mu}$) are called {\it dual} if they satisfy the following properties:
    \begin{itemize}
        \item[(i)] skew Cauchy identity
        $$
       \sum_{\lambda} F_{\lambda/\mu}(\mathbf{x}) f_{\lambda/\nu}(\mathbf{y}) =  \Omega(\mathbf{x},\mathbf{y}) \sum_{\kappa} F_{\nu/\kappa}(\mathbf{x}) f_{\mu/\kappa}(\mathbf{y}),
        $$
        for some {\it Cauchy kernel} $\Omega %= \prod_{i,j} P(x_i y_j) 
        %\in \hat\Lambda \otimes \hat\Lambda
        $ 
        %($P(t)$ is some formal power series in $t$) 
        such that there is an automorphism $\widetilde{\omega}: \hat\Lambda \to \hat\Lambda$ satisfying $\widetilde{\omega}_{\mathbf{x}} : \Omega(\mathbf{x}, \mathbf{y}) \mapsto \Omega(\mathbf{x}, \mathbf{y})^{-1}$
        \item[(ii)] the elements $\{f_{\lambda}\}$ are linearly independent, where $f_{\lambda} := f_{\lambda/\varnothing}$. 
        %\item[(iii)]
    \end{itemize}
\end{definition}

\begin{lemma}[Skew Pieri-type formula]\label{sptf}
    Let $\{ F_{\lambda/\mu}\}, \{f_{\lambda/\mu} \}$ be dual families of symmetric functions. Then the following formula holds:
   \begin{align}
            \sum_{\lambda, \eta} \widetilde{F}_{\nu/\eta}(\mathbf{x}) F_{\lambda/\mu}(\mathbf{x}) f_{\lambda/\eta}(\mathbf{y}) &= \Omega(\mathbf{x}, \mathbf{y}) f_{\mu/\nu}(\mathbf{y}), \qquad \widetilde{F}_{\nu/\eta}(\mathbf{x}) := \widetilde{\omega}(F_{\nu/\eta}(\mathbf{x})).
    \end{align}
\end{lemma}
\begin{proof}
Specializing $\nu = \varnothing$ in the skew Cauchy identity we obtain the dual Pieri-type formulas
$$
        \sum_{\lambda} F_{\lambda/\mu}(\mathbf{x}) f_{\lambda}(\mathbf{y}) =  \Omega(\mathbf{x},\mathbf{y}) f_{\mu}(\mathbf{y}), \qquad \sum_{\lambda} \widetilde{F}_{\lambda/\mu}(\mathbf{x}) f_{\lambda}(\mathbf{y}) =  \Omega(\mathbf{x},\mathbf{y})^{-1} f_{\mu}(\mathbf{y}).
$$
Using these formulas we have
\begin{align*}
    \sum_{\lambda, \rho} F_{\lambda/\rho}(\mathbf{x}) \widetilde{F}_{\rho/\mu}(\mathbf{x}) f_{\lambda}(\mathbf{y}) &= \Omega(\mathbf{x}, \mathbf{y}) \sum_{\rho}  \widetilde{F}_{\rho/\mu}(\mathbf{x}) f_{\rho}(\mathbf{y})  
    = \Omega(\mathbf{x}, \mathbf{y}) \Omega(\mathbf{x}, \mathbf{y})^{-1} \ f_{\mu}(\mathbf{y}) 
    = f_{\mu}(\mathbf{y}).
\end{align*}
Since the family $\{ f_{\lambda}\}$ is linearly independent we conclude the following orthogonality relation
$$
    \sum_{\rho} F_{\lambda/\rho}(\mathbf{x}) \widetilde{F}_{\rho/\mu}(\mathbf{x}) = \delta_{\lambda \mu}. \eqno{(o)}
$$
Now we have the following
\begin{align*}
    \sum_{\lambda, \eta} \widetilde{F}_{\nu/\eta}(\mathbf{x}) F_{\lambda/\mu}(\mathbf{x}) f_{\lambda/\eta}(\mathbf{y}) &= \sum_{\eta} \widetilde{F}_{\nu/\eta}(\mathbf{x}) \sum_{\lambda} F_{\lambda/\mu}(\mathbf{x}) f_{\lambda/\eta}(\mathbf{y}) \\
    &= \Omega(\mathbf{x}, \mathbf{y}) \sum_{\eta} \widetilde{F}_{\nu/\eta}(\mathbf{x}) \sum_{\kappa} F_{\eta/\kappa}(\mathbf{x}) f_{\mu/\kappa}(\mathbf{y})\\
    &= \Omega(\mathbf{x}, \mathbf{y}) \sum_{\kappa} f_{\mu/\kappa}(\mathbf{y}) \underbrace{\sum_{\eta}\widetilde{F}_{\nu/\eta}(\mathbf{x}) F_{\eta/\kappa}(\mathbf{x})}_{=\delta_{\nu \kappa} \text{ by ($o$)}} \\
    &= \Omega(\mathbf{x}, \mathbf{y}) f_{\mu/\nu}(\mathbf{y}),
\end{align*}
%where at the last step we used the orthogonality relation just derived above. 
\end{proof}
\begin{remark}
This is a general formulation of the method in  \cite{warnaar} for Hall-Littlewood polynomials.
\end{remark}
\subsection{Back to Grothendieck} 
Recall that $J^{\beta}_{\lambda/\!\!/\mu} = \omega(G^{\beta}_{\lambda/\!\!/\mu})$ and $j^{\beta}_{\lambda/\mu} = \omega(g^{\beta}_{\lambda/\mu})$, where $\omega$ is the standard involution defined on the Schur basis by $\omega : s_{\lambda} \mapsto s_{\lambda'}$ (in case of $G \in \hat\Lambda$ it is extended for infinite linear combinations). 
\begin{theorem}[Skew Pieri-type formulas] We have
\begin{align*}
    \sum_{\lambda, \eta}   g^{\beta}_{\lambda/\mu}(\mathbf{y}) j^{\beta}_{\nu/\eta}(-\mathbf{y}) G^{\beta}_{\lambda/\!\!/\eta}(\mathbf{x}) &= \prod_{i,j} \frac{1}{(1 - x_i y_j)} G^{\beta}_{\mu/\!\!/\nu}(\mathbf{x}), \\ 
    \sum_{\lambda, \eta}  G^{\beta}_{\lambda/\!\!/\mu}(\mathbf{x}) J^{\beta}_{\nu/\!\!/\eta}(-\mathbf{x}) g^{\beta}_{\lambda/\eta}(\mathbf{y}) &= \prod_{i,j} \frac{1}{(1 - x_i y_j)} g^{\beta}_{\mu/\nu}(\mathbf{y})
\end{align*}
and the following dual formulas
\begin{align*}
   \sum_{\lambda, \eta}   j^{\beta}_{\lambda/\mu}(\mathbf{y}) g^{\beta}_{\nu/\eta}(-\mathbf{y}) G^{\beta}_{\lambda/\!\!/\eta}(\mathbf{x}) &= \prod_{i,j} {(1 + x_i y_j)} G^{\beta}_{\mu/\!\!/\nu}(\mathbf{x}),\\
    \sum_{\lambda, \eta} J^{\beta}_{\lambda/\!\!/\mu}(\mathbf{x}) G^{\beta}_{\nu/\!\!/\eta}(-\mathbf{x})  g^{\beta}_{\lambda/\eta}(\mathbf{y}) &= \prod_{i,j} {(1 + x_i y_j)} g^{\beta}_{\mu/\nu}(\mathbf{y}).
\end{align*}
\end{theorem}
\begin{proof}
The formulas imply from Lemma \ref{sptf} for corresponding dual skew families $G, g, J, j$ and the automorphism $\widetilde{\omega} : f(\mathbf{x}) \mapsto \omega(f)(-\mathbf{x})$ satisfying $\widetilde{\omega} \prod_{i,j} (1 - x_i y_j)^{-1} = \prod_{i,j} (1 - x_i y_j)$. 
\end{proof}
\begin{corollary}
Specializing for single variables we have 
\begin{align}
    \prod_{i} \frac{1}{(1 - x_i y)} G^{\beta}_{\mu/\!\!/\nu}(\mathbf{x}) &= \sum_{{\lambda \supset \mu \atop \nu/\eta \text{ vert. strip}}}   g^{\beta}_{\lambda/\mu}(y) j^{\beta}_{\nu/\eta}(-y) G^{\beta}_{\lambda/\!\!/\eta}(\mathbf{x}) \nonumber\\ &=\sum_{{\lambda \supset \mu \atop \nu/\eta \text{ vert. strip}}}  y^{c(\lambda/\mu) + c(\nu/\eta)}\beta^{|\lambda/\mu| - c(\lambda/\mu)}(\beta - y)^{|\nu/\eta| - c(\nu/\eta)} G^{\beta}_{\lambda/\!\!/\eta}(\mathbf{x}),\label{xyg} \\ 
    \prod_{j} \frac{1}{(1 - x y_j)} g^{\beta}_{\mu/\nu}(\mathbf{y}) &= 
    \sum_{{\lambda/\mu \text{ hor strip} \atop \nu/\eta \text{ vert strip}}}  G^{\beta}_{\lambda/\!\!/\mu}({x}) J^{\beta}_{\nu/\!\!/\eta}(-{x}) g^{\beta}_{\lambda/\eta}(\mathbf{y}), \nonumber \\
    &= \sum_{{\lambda/\mu \text{ hor strip} \atop \nu/\eta \text{ vert strip}}}   (-1)^{|\nu/\eta|} (1 - \beta x)^{a(\lambda/\!\!/\mu) - a(\nu'/\!\!/\eta') - |\nu/\eta|} x^{|\lambda/\mu| + |\nu/\eta|} g^{\beta}_{\lambda/\eta}(\mathbf{y}),
\end{align}
as well as the dual formulas
\begin{align}
    \prod_{i} {(1 + x_i y)} G^{\beta}_{\mu/\!\!/\nu}(\mathbf{x}) &= \sum_{{\lambda/\mu \text{ vert strip} \atop \eta \subset \nu}}   j^{\beta}_{\lambda/\mu}({y}) g^{\beta}_{\nu/\eta}(-{y}) G^{\beta}_{\lambda/\!\!/\eta}(\mathbf{x}),\nonumber \\
    &= \sum_{{\lambda/\mu \text{ vert strip} \atop \eta \subset \nu}}  (-1)^{c(\nu/\eta)} y^{c(\lambda/\mu) + c(\nu/\eta)} (\beta + y)^{|\lambda/\mu| - c(\lambda/\mu)} \beta^{|\nu/\eta| - c(\nu/\eta)} G^{\beta}_{\lambda/\!\!/\eta}(\mathbf{x}),\\
\prod_{j} {(1 + x y_j)} g^{\beta}_{\mu/\nu}(\mathbf{y}) &=
    \sum_{{\lambda/\mu \text{ vert strip} \atop \nu/\eta \text{ hor strip}}}  J^{\beta}_{\lambda/\!\!/\mu}({x}) G^{\beta}_{\nu/\!\!/\eta}(-{x}) g^{\beta}_{\lambda/\eta}(\mathbf{y}),\nonumber\\
    &=\sum_{{\lambda/\mu \text{ vert strip} \atop \nu/\eta \text{ hor strip}}}  (-1)^{|\nu/\eta|}(1 + \beta x)^{-a(\lambda'/\!\!/\mu') - |\lambda/\mu| + a(\nu/\!\!/\eta)} x^{|\lambda/\mu| + |\nu/\eta|}  g^{\beta}_{\lambda/\eta}(\mathbf{y}).
\end{align}
Hence we also obtain the following skew rules
\begin{align}
\label{hg} h_k G^{\beta}_{\mu/\!\!/\nu} &= \sum_{{\lambda \supset \mu \atop \nu/\eta \text{ vert. strip}}}  (-1)^{k - c(\lambda/\mu) - c(\nu/\eta)} \beta^{|\lambda/\mu| - k} \binom{|\nu/\eta| - c(\nu/\eta)}{k -c(\lambda/\mu) - c(\nu/\eta)} G^{\beta}_{\lambda/\!\!/\eta},\\
\label{hgg} h_k g^{\beta}_{\mu/\nu} &= \sum_{{\lambda/\mu \text{ hor strip} \atop \nu/\eta \text{ vert strip}}}  (-1)^{k - |\lambda/\mu|} \beta^{k - |\lambda/\mu| - |\nu/\eta|} \binom{a(\lambda/\!\!/\mu) - a(\nu'/\!\!/\eta') - |\nu/\eta|}{k - |\lambda/\mu| - |\nu/\eta|} g^{\beta}_{\lambda/\eta},\\
\label{eg} e_k G^{\beta}_{\mu/\!\!/\nu} &= \sum_{{\lambda/\mu \text{ vert strip} \atop \eta \subset \nu}}  (-1)^{c(\nu/\eta)} \beta^{|\lambda/\mu| + |\nu/\eta| - k} \binom{|\lambda/\mu| - c(\lambda/\mu)}{k - c(\lambda/\mu) - c(\nu/\eta)} G^{\beta}_{\lambda/\!\!/\eta},\\
\label{egg} e_k g^{\beta}_{\mu/\nu} &= \sum_{{\lambda/\mu \text{ vert strip} \atop \nu/\eta \text{ hor strip}}}  (-1)^{|\nu/\eta|} \beta^{k - |\lambda/\mu| - |\nu/\eta|} \binom{-a(\lambda'/\!\!/\mu') - |\lambda/\mu| + a(\nu/\!\!/\eta)}{k - |\lambda/\mu| - |\nu/\eta|} g^{\beta}_{\lambda/\eta}.
\end{align}
\end{corollary}

\begin{proof}[Proof of Theorem \ref{spieri}]
Using the Cauchy identity and the fact that $j^{\beta}_{\rho}(y) = 0$ unless $\rho$ is a single column and $j^{\beta}_{(1^k)}(y) = y(\beta+y)^{k-1}$, we obtain
$$
\sum_{\rho} G^{\beta}_{\rho}(\mathbf{x}) j^{\beta}_{\rho}(y) = 1+ \sum_{k \ge 1} G^{\beta}_{(1^k)}(\mathbf{x}) y(\beta+y)^{k-1} = \prod_{i} (1 + x_i y).
$$
Therefore using \eqref{xyg} we have
\begin{align*}
 \prod_{i} (1 + x_i y) G^{\beta}_{\mu/\!\!/\nu}(\mathbf{x}) &= \left( 1+ \sum_{k \ge 1} G^{\beta}_{(1^k)}(\mathbf{x}) y(\beta+y)^{k-1}\right)G^{\beta}_{\mu/\!\!/\nu}(\mathbf{x})\\
 &= \sum_{{\lambda/\mu \text{ vert strip} \atop \eta \subset \nu}}  (-1)^{c(\nu/\eta)} y^{c(\lambda/\mu) + c(\nu/\eta)} (\beta + y)^{|\lambda/\mu| - c(\lambda/\mu)} \beta^{|\nu/\eta| - c(\nu/\eta)} G^{\beta}_{\lambda/\!\!/\eta}(\mathbf{x})
\end{align*}
The expansion for $G^{\beta}_{(1^k)} G^{\beta}_{\mu/\!\!/\nu}$ is obtained by comparing the coefficients at $[y(\beta + y)^{k-1}]$ from both sides.
Notice that $g^{\beta}_{(k)} = h_k$ and hence \eqref{hgg} gives the needed skew Pieri rule for $g^{\beta}_{(k)} g^{\beta}_{\mu/\nu}$.
\end{proof}

\section{Basis phenomenon}\label{bf}
A natural question is when a family of symmetric functions of unbounded degree, i.e., belonging to the completion of the ring of symmetric functions forms a {\it basis} of a certain ring, like $\{ G_{\lambda}\}$ does. In this section we prove that products $G_{\mu} G_{\nu}$ expand finitely in $\{ G_{\lambda}\}$ {without} appealing to a Littlewood-Richardson rule as in \cite{buch, pp1}. We provide a general sufficient condition for this situation.

Say that a family $\{ f_{\lambda}\}$ %(of elements of a completion of a certain algebra) 
presents a {\it basis phenomenon} if its elements are linearly independent and for all $\mu,\nu$, the product expansion $f_{\mu} f_{\nu} = \sum_{\lambda} c^{\lambda}_{\mu \nu} f_{\lambda}$ exists and is {\it finite}. 

\begin{theorem}\label{gfin}
The family $\{G_{\lambda}\}$ presents a basis phenomenon. %, i.e., %of the commutative ring ${\Gamma} = \bigoplus_{\lambda} \mathbb{Z} \cdot G_{\lambda}$. 
%the product expansion 
%$$
%G_{\mu} G_{\nu} = \sum_{\lambda} c^{\lambda}_{\mu \nu} G_{\lambda}, \qquad |\lambda| \ge |\mu| + |\nu|
%$$
%is finite.
\end{theorem}

%\subsection{A sufficient condition}
\begin{definition}
Say that a family $\{ F_{\lambda} \}$ %, $\{ f_{\lambda}\}$ 
of symmetric functions $F_{\lambda}\in \hat\Lambda$ (assume $F_{\varnothing} = 1$ for simplicity) is {\it damping} if % $f_{\lambda} \in \Lambda$, {\it dual} via a certain Cauchy identity
%$$
%\sum_{\lambda} F_{\lambda}(\mathbf{x}) f_{\lambda}(\mathbf{y}) = \Omega(\mathbf{x}, \mathbf{y}).
%$$ 
\begin{itemize}
\item[(i)] it is linearly independent in $\hat\Lambda$ and each element of $\hat\Lambda$ can uniquely be expressed as (possibly) an infinite linear combination of $\{F_{\lambda} \}$
\item[(ii)] there is an involutive automorphism $\widehat{\omega} : F_{\lambda} \mapsto F_{\lambda'}$
\item[(iii)] there is a Pieri-type formula\footnote{Which is a special case of a skew Cauchy for $\mu = \varnothing$.}
\begin{align}\label{dpieri}
\sum_{\lambda} F_{\lambda}(\mathbf{x}) f_{\lambda \mu}(\mathbf{y}) = \Omega(\mathbf{x}, \mathbf{y}) F_{\mu}(\mathbf{x})
\end{align}
for some (non-degenerate) kernel $\Omega$ and a {\it dual} family $\{f_{\lambda \mu}\}$ that satisfies % \in \hat{\Lambda} \otimes \hat{\Lambda}
\begin{itemize}
\item[(a)] the {\it damping condition}: 
For each $\mu$ and $n \in \mathbb{N}$, there exists a constant $k = k(\mu, n)$ such that if $f_{\lambda \mu}(x_1, \ldots, x_n) \ne 0$ then $\lambda_1 < k$ 
\item[(b)] $\{f_{\lambda} \}$ are linearly independent (where $f_{\lambda} := f_{\lambda \varnothing}$). 
\end{itemize}
\end{itemize}
\end{definition}
\begin{remark}
The damping condition is important here. It is a natural property of symmetric functions whose combinatorial presentations have {\it strict} row conditions. For example, Schur polynomials satisfy it: $s_{\lambda'/\mu'}(x_1, \ldots, x_n) \ne 0$ implies that $\lambda_1 \le k(\mu, n) = \mu_1 + n$. On the other hand, the dual Grothendieck polynomials $g_{\lambda/\mu}$ do {\it not} have the damping condition, i.e., if $g_{\lambda/\mu}(x_1, \ldots, x_n) \ne 0$ then $\lambda_1$ (and $\ell(\lambda)$) can be arbitrarily large regardless of $\mu$ and $n$, as their formula is based on plane partitions that have weak inequalities in rows and columns. 
\end{remark}
\begin{theorem}
Every damping family of symmetric functions presents a basis phenomenon.
\end{theorem}
\begin{proof}
Let $\{F_{\lambda} \}$ be a damping family. By \eqref{dpieri}, combining it for $\mu = \varnothing$, we have
\begin{align*}
\sum_{\lambda} F_{\lambda}(\mathbf{x}) f_{\lambda \mu}(\mathbf{y}) = \Omega(\mathbf{x}, \mathbf{y}) F_{\mu}(\mathbf{x}) = \sum_{\nu} F_{\nu}(\mathbf{x}) f_{\nu}(\mathbf{y}) F_{\mu}(\mathbf{x}).
\end{align*}
Therefore,
$$
F_{\mu} F_{\nu} = \sum_{\lambda} d^{\lambda}_{\mu \nu} F_{\lambda} \implies f_{\lambda \mu} = \sum_{\nu} d^{\lambda}_{\mu \nu} f_{\nu}.
$$
For any fixed $\mu, \nu$ consider $\lambda$ so that $d^{\lambda}_{\mu \nu} \ne 0$. Let $n = n(\nu)$ be a minimal number such that $f_{\nu}(x_1, \ldots, x_n) \ne 0$. By the damping condition we have 
$$
f_{\lambda \mu}(x_1, \ldots, x_n) = \sum_{\nu} d^{\lambda}_{\mu \nu} f_{\nu}(x_1, \ldots, x_n) \ne 0 \implies \lambda_1 < k = k(\mu, n(\nu)),
$$
which means that  $\lambda_1$ is bounded from above by a constant depending on $\mu, \nu$. Since there is an automorphism $\widehat{\omega}$ mapping $F_{\lambda}$ to $F_{\lambda'},$ we also have $d^{\lambda}_{\mu \nu} = d^{\lambda'}_{\mu' \nu'} \ne 0$ and hence by the same argument we obtain that $\lambda'_1 = \ell(\lambda) < k(\mu',n(\nu'))$ is bounded from above as well.  Therefore, for every pair $\mu, \nu$ there are only {\it finitely} many $\lambda$ so that $d^{\lambda}_{\mu \nu} \ne 0$. 
\end{proof}

Consider a damping family $\{F_{\lambda}\}$ with a dual $\{f_{\lambda \mu} \}$, so that there is an automorphism $f_{\lambda} \mapsto f_{\lambda'}$. Define the ring $\Phi := \bigoplus_{\lambda} K \cdot F_{\lambda}$ with the basis $\{ F\}$. Suppose $\{F\}$ has a skew extension $\{F_{\lambda \mu} \}$ so that the following dual formula holds as well: 
$$
\sum_{\lambda} F_{\lambda \mu}(\mathbf{x}) f_{\lambda}(\mathbf{y}) = \Omega(\mathbf{x}, \mathbf{y}) f_{\mu}(\mathbf{y}).
$$
\begin{proposition}\label{phi}
$F_{\lambda \mu} \in \Phi$.
\end{proposition}
\begin{proof}
The proof is similar to the previous result. The dual formula gives 
$$F_{\lambda \mu} =  \sum_{\nu} d^{\mu \nu}_{\lambda} F_{\nu} \implies f_{\mu} f_{\nu} = \sum_{\lambda} d^{\mu \nu}_{\lambda} f_{\lambda}.$$ 
For fixed $\lambda, \mu$ consider $\nu$ so that $d^{\mu \nu}_{\lambda} \ne 0$ and let $n = n(\lambda)$ be a minimal number so that $f_{\lambda}(x_1, \ldots, x_n) \ne 0$. Then $f_{\mu}(x_1, \ldots, x_n) f_{\nu}(x_1, \ldots, x_n) \ne 0$ and hence from the damping condition we have $\nu_1 < k = k(n(\lambda))$ is bounded. Since $d^{\mu \nu}_{\lambda} = d^{\mu' \nu'}_{\lambda'}$ by the same argument we have $\nu'_1 = \ell(\nu) < k(n(\lambda'))$ is bounded from above as well. Therefore, for every $\lambda, \mu$ there exists only finitely many $\nu$ for which $d^{\mu \nu}_{\lambda} \ne 0$.
\end{proof}

\begin{proof}[Proof of Theorem \ref{gfin}]
Let us show that $\{G_{\lambda} \}$ is a damping family. The conditions (i), (ii) are satisfied (see e.g. \cite{buch, dy}). For (iii) we take the following Pieri-type formula that we obtained earlier\footnote{We did not take another dual formula containing $g_{\lambda}$, since it does not satisfy the damping condition.}
$$
\sum_{\lambda} G_{\lambda}(\mathbf{x}) j_{\lambda/\mu}(\mathbf{y}) = \prod_{i,j} (1 + x_i y_j) G_{\mu}(\mathbf{x}).
$$
The family $\{j_{\lambda/\mu}\}$ satisfies the damping condition: by definition (see Sec. \ref{jjj}) $j_{\lambda/\mu}$ has a combinatorial formula over certain tableaux that are {\it row strict} (the operator $\overline{B}(x) \lambda$ removes vertical strips from $\lambda$), which means that if $j_{\lambda/\mu}(x_1, \ldots, x_n) \ne 0$ then $\lambda_1 \le \mu_1 + n$ is bounded from above. 
\end{proof}

\begin{corollary}
Let $\Gamma := \bigoplus_{\lambda} \mathbb{Z} \cdot G_{\lambda}$. We have $G_{\lambda/\!\!/\mu} \in \Gamma$ and $G_{\lambda/\mu} \in \Gamma$.
\end{corollary}
\begin{proof}
For $G_{\lambda/\!\!/\mu}$, the result follows by Proposition \ref{phi} since $\{G_{\lambda}\}$ is damping and for $G_{\lambda/\mu}$ since $G_{\lambda/\mu} = \sum_{\nu \subset \mu} G_{\lambda/\!\!/\nu}$ is a finite sum. 
\end{proof}

%Another example: Consider the functions $S^{\alpha}_{\lambda}(x_1, x_2,\ldots) := s_{\lambda}(x_1/(1 - \alpha x_1), x_2/(1 - \alpha x_2), \ldots)$. These functions appear as a special case of a canonical two parameter deformation of Grothendieck polynomials $G^{(\alpha, \beta)}_{\lambda}$ \cite{dy}. It is clear that $S^{\alpha}_{\lambda}$ multiplies via the usual Littlewood-Richardson coefficients, but let us show that it is an S-basis. We have the following Pieri-type formula
%$$
%\sum_{\lambda} S^{\alpha}_{\lambda}(\mathbf{x}) s_{\lambda'/\mu'}(\mathbf{y}) = \prod_{i,j}\left(1 + \frac{x_i}{1 - \alpha x_i} y_j \right) s_{\mu'}(\mathbf{y}).
%$$
%Clearly, the function $\widetilde{s}_{\lambda/\mu} = s_{\lambda'/\mu'}$ satisfies non-vanishing bound property.
%$\{S_{\lambda} \}$

\section{Dual filtered Young graphs}
Following \cite{pp}, a {\it weighted filtered graph} is a digraph $G = (V, r, E, w)$ where $V$ is a set of countably many vertices together with a {\it rank function} $r : V \to \mathbb{Z}$ satisfying $r(a) \le r(b)$ for every (directed) edge $(a, b) \in E$, and $w : E \to \mathbb{R}$ is some weight function. %Similarly, a {\it strict filtered graph} satisfies $r(a) < r(b)$. 

For a pair $G_1 = (V, r, E_1, w_1), G_2 = (V, r, E_2, w_2)$ of filtered graphs on the same (ranked) vertex set $V$ construct a digraph $G  = (V, E)$ so that $E = E_1 \cup \overline{E}_2$ is a union of edges $E_1$ and edges of $E_2$ but taken in {\it opposite} direction. Let $\mathbb{R}V$ be the free abelian group on $V$ (formal $\mathbb{R}$-linear combinations of vertices $V$). Define the {\it up} and {\it down} operators $U, D \in End(\mathbb{R}V)$ on $G$ as follows:
$$
U v = \sum_{e = (v \to u) \in E} w_1(e) u, \qquad\qquad D v = \sum_{e = (u \to v) \in E} w_2(e) u.
$$
Say that $G$ is a {\it dual filtered graph} if there exist scalars $\alpha, \beta \in \mathbb{R}$ such that for all $v \in V$ we have
$$
[D,U] v = (DU - UD)  v = (\alpha + \beta D) v.
$$
\begin{remark}
Up to normalizations, there are only three distinct types of $(\alpha, \beta) \in \{(1,1), (0,1), (1,0)\}$. In addition, the relation $[D, U] = D$ can be shifted with $D' = D - 1$ which gives $[D', U] = 1 + D'$.
\end{remark}
\begin{remark}
If the rank function $r$ satisfies $r(a) + 1 = r(b)$ for every edge $(a,b)$ and $(\alpha, \beta) = (1,0)$, then the corresponding graphs $G$ are called {\it dual graded graphs} studied by Fomin \cite{fomindual} and by Stanley \cite{stadiff} as {\it differential posets}.
\end{remark}
\begin{remark}
%Let $K$ be a field. 
An associative algebra generated by $U, D$ subject to $[D, U] = 1$ is called the {\it first Weyl algebra}. One may consider $D= \frac{d}{dx}$ as a differential operator and $U = x$ acting on a polynomial ring $K[x]$. The relation $[D, U] = 1 + D$ corresponds to the difference operator $D f(x) = f(x+1) - f(x)$.
\end{remark}

\subsection{New constructions of filtered Young graphs} 
Recall that $\mathbb{Y}$ is the Young lattice, i.e., an infinite graph whose vertices are indexed by partitions and edges are given by $(\lambda, \lambda + \square)$. We think of $\mathbb{Y}$ as a self-dual graph with up and down directed edges $(\lambda \to \lambda \pm \square)$. 

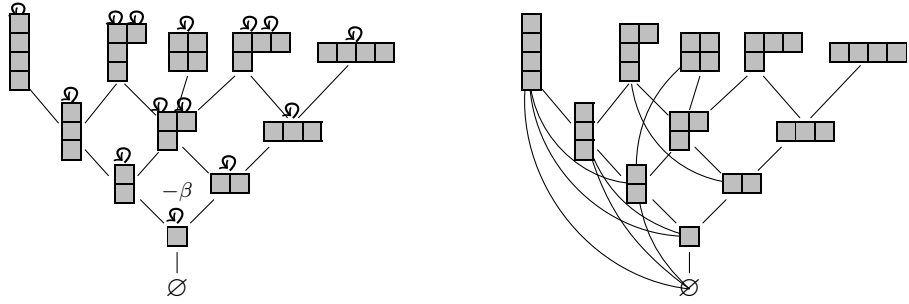
\begin{figure}[t]
\begin{center}
\ytableausetup{smalltableaux}
{\scriptsize
\begin{tikzpicture}[scale = 0.7]
\draw (-2.25, 2.15) to (-2.80, 2.8); \draw (-1.75, 2.15) to (-1.15, 2.9); \draw (-0.40, 2.30) to (-1.0, 2.9); \draw (-0.0, 2.40) to (0.2, 3.05); \draw (0.4, 2.4) to (1.1, 3.05); \draw (2.1, 2.25) to (1.4, 3.05); \draw (2.3, 2.25) to (3.25, 3.25);

\draw (-1.30, 1.15) to (-1.75, 1.65); \draw (-0.75, 1.15) to (-0.35, 1.60); \draw (0.60, 1.20) to (0.10, 1.70); \draw (1.3, 1.25) to (1.75, 1.70);

\draw (-0.25, 0.25) to (-0.7, 0.7); \draw (0.25, 0.25) to (0.7, 0.7);

\draw (0,-0.7) to (0,-0.3);

%	\vertex (1) at (0,2) [label=above left:$1$]{};
%	\vertex (2) at (2.5,2) [label=above right:$2$]{};	

%\draw[->] (1) edge[bend left=20] node[above]{$e_1$} (2);

\draw (0,0.25)  edge[loop, thick] node[above]{\scriptsize$-\beta$} (0,0.25);
\draw (-1,1.40)  edge[loop, thick] node[above]{} (-1,1.40);
\draw (1,1.25)  edge[loop, thick] node[above]{} (1,1.25);
\draw (-2,2.55)  edge[loop, thick] node[above]{} (-2,2.55);

\draw (-0.3,2.35)  edge[loop, thick] node[above]{} (-0.3,2.35); %
\draw (0.15,2.35)  edge[loop, thick] node[above]{} (0.15,2.35);

\draw (2.2,2.25)  edge[loop, thick] node[above]{} (2.2,2.25);
\draw (-3,4.15)  edge[loop, thick] node[above]{} (-3,4.15);

\draw (-1.15,4.0)  edge[loop, thick] node[above]{} (-1.15,4.0); %\scriptsize$-2$
\draw (-0.75,4.0)  edge[loop, thick] node[above]{} (-0.75,4.0);

\draw (0.2,3.9)  edge[loop, thick] node[above]{} (0.2,3.9);

\draw (1.4,3.9)  edge[loop, thick] node[above]{} (1.4,3.9);
\draw (1.8,3.9)  edge[loop, thick] node[above]{} (1.8,3.9);

\draw (3.4,3.7)  edge[loop, thick] node[above]{} (3.4,3.7);

%\draw (0,0.5)  edge[bend right=100] (0,0.1);

\node at (-3,3.5) {\ydiagram[*(lightgray)]{1,1,1,1}}; \node at (-0.95,3.5) {\ydiagram[*(lightgray)]{2,1,1}}; \node at (0.2,3.5) {\ydiagram[*(lightgray)]{2,2}}; \node at (1.6,3.5) {\ydiagram[*(lightgray)]{3,1}}; \node at (3.4,3.5) {\ydiagram[*(lightgray)]{4}}; 
\node at (-2,2) {\ydiagram[*(lightgray)]{1,1,1}}; \node at (0,2) {\ydiagram[*(lightgray)]{2,1}}; \node at (2.2,2) {\ydiagram[*(lightgray)]{3}};
\node at (-1,1) {\ydiagram[*(lightgray)]{1,1}}; \node at (1,1) {\ydiagram[*(lightgray)]{2}};
\node at (0,0) {\ydiagram[*(lightgray)]{1}};
\node at (0,-1) {\normalsize$\varnothing$};
\end{tikzpicture}
\qquad\qquad 
\begin{tikzpicture}[scale = 0.7]
\draw (-2.25, 2.15) to (-2.80, 2.8); \draw (-1.75, 2.15) to (-1.15, 2.9); \draw (-0.40, 2.30) to (-1.0, 2.9); \draw (-0.0, 2.40) to (0.2, 3.05); \draw (0.4, 2.4) to (1.1, 3.05); \draw (2.1, 2.25) to (1.4, 3.05); \draw (2.3, 2.25) to (3.25, 3.25);

\draw (-1.30, 1.15) to (-1.75, 1.65); \draw (-0.75, 1.15) to (-0.35, 1.60); \draw (0.60, 1.20) to (0.10, 1.70); \draw (1.3, 1.25) to (1.75, 1.70);

\draw (-0.25, 0.25) to (-0.7, 0.7); \draw (0.25, 0.25) to (0.7, 0.7);

\draw (0,-0.7) to (0,-0.3);

\draw (-1,1)  edge[bend right=20] (0,-1);
\draw (-2,2)  edge[bend right=20] (0,-1); \draw (-2,2)  edge[bend right=30] (0,0); 
\draw (-3,3.5)  edge[bend right=50] (0,-1); \draw (-3,3.5)  edge[bend right=50] (0,0); \draw (-3,3.5)  edge[bend right=50] (-1,1); 
\draw (-1.1,2.9)  edge[bend right=40] (1,1); 
\draw (0.2,3.5)  edge[bend right=30] (-1,1); 

\node at (-3,3.5) {\ydiagram[*(lightgray)]{1,1,1,1}}; \node at (-0.95,3.5) {\ydiagram[*(lightgray)]{2,1,1}}; \node at (0.2,3.5) {\ydiagram[*(lightgray)]{2,2}}; \node at (1.6,3.5) {\ydiagram[*(lightgray)]{3,1}}; \node at (3.4,3.5) {\ydiagram[*(lightgray)]{4}}; 

\node at (-2,2) {\ydiagram[*(lightgray)]{1,1,1}}; \node at (0,2) {\ydiagram[*(lightgray)]{2,1}}; \node at (2.2,2) {\ydiagram[*(lightgray)]{3}};

\node at (-1,1) {\ydiagram[*(lightgray)]{1,1}}; \node at (1,1) {\ydiagram[*(lightgray)]{2}};

\node at (0,0) {\ydiagram[*(lightgray)]{1}};

\node at (0,-1) {\normalsize$\varnothing$};

\end{tikzpicture}

}
\end{center}
\caption{The dual filtered Young graph $\beta{\mathbb{Y}}$. The graph on the left corresponds to {\it up} edges and on the right to {\it down} edges. Here each loop has the weight $-\beta$. }\label{yg}
\end{figure}

{\bf I. } First define the following {\it $\beta$-filtration} $\beta{\mathbb{Y}}$ of Young's lattice $\mathbb{Y}$ (see Fig. \ref{yg}):
\begin{itemize}
\item[(i)] vertices $V$ are integer partitions ranked by the number of boxes $r(\lambda) = |\lambda|$ 
\item[(ii)] {\it up} edges (of $E_1$) are as in Young's lattice $(\lambda \to \lambda + \square)$ with the weight $w = 1$ but there are also $i(\lambda)$ many {\it loops} $(\lambda \to \lambda)$ each with the weight $w = -\beta$ (recall that $i(\lambda)$ is the number of inner corners of $\lambda$)
\item[(iii)] {\it down} edges (of $\overline{E}_2$) are given by $(\lambda \to \mu)$ iff all boxes $\lambda/\mu$ are on a single column, and the corresponding weight is $w = \beta^{|\lambda/\mu| - 1}.$
\end{itemize}

{\bf II. } Next, let $\varkappa$ be a scalar parameter and define the {\it Cauchy filtration} of ${\mathbb{Y}}$ denoted by $\varkappa{\mathbb{Y}}$ that satisfies exactly the same conditions (i) and (ii) as $\beta{\mathbb{Y}}$ but its {\it down} edges are given by $(\lambda \to \mu)$ iff $\lambda \supset \mu$ with the weight $w = \varkappa^{c(\lambda/\mu)} \beta^{|\lambda/\mu| - c(\lambda/\mu)}$.

\begin{theorem} We have
\begin{itemize}
\item[(i)] $\beta{\mathbb{Y}}$ is a dual filtered graph satisfying $[D, U] = 1$.
\item[(ii)] $\varkappa{\mathbb{Y}}$ is a dual filtered graph satisfying $[D, U] = \varkappa(1 + D)$.
\end{itemize}
\end{theorem}
\begin{proof}
These constructions are natural consequences of the Cauchy identity. 
Suppose the operator series
$$A(x) = 1 + \sum_{i  \ge 1} U_i x^i, \qquad B(y) = 1 + \sum_{i \ge 1} D_i y^i,$$
satisfy the Cauchy identity
$$
B(y) A(x) = (1 - xy)^{-1} A(x) B(y). % \frac{1}{1 - xy}
$$
By comparing coefficients at $xy$ and $x$ after plugging $y = \varkappa$ we obtain that
$$[D_1, U_1] = 1 \qquad\text{ and }\qquad [{D}(\varkappa), U_1] = \varkappa(1 + {D}(\varkappa)),\quad {D}(\varkappa) = D_1 \varkappa + D_2 \varkappa^2 + \cdots.$$
Let $\beta \in \mathbb{R}$ and define the operators 
\begin{align}\label{opers}
\widetilde{U} = \widetilde{u}_1 + \widetilde{u}_2 + \cdots, \qquad \widetilde{D} = \widetilde{d}_1 + \widetilde{d}_2 + \cdots, \qquad \overline{D} = -1 + (1 + \varkappa\widetilde{d}_1)(1 + \varkappa\widetilde{d}_2) \cdots 
\end{align}
Then the Cauchy identity gives 
$
[\widetilde{D},\widetilde{U}] = 1
$
and $[\overline{D}, \widetilde{U}] = \varkappa(1 + \overline{D})$. Observe that $\widetilde{U}$ defines the up edges of $\beta\mathbb{Y}$ and $\varkappa\mathbb{Y}$, the operator $\widetilde{D}$ corresponds to the down operator of $\beta\mathbb{Y}$ as it defines its down edges. The operator $\overline{D}$ defines the down edges of $\varkappa\mathbb{Y}$ as $\varkappa^k \widetilde{d}_{i_1} \cdots \widetilde{d}_{i_k}$ %for any $i_1 < \cdots < i_k$ 
removes boxes from the $k$ columns $i_1 < \cdots < i_k$ in all possible ways giving the corresponding weight $\varkappa^k \beta^{|\lambda/\mu| - k}$. 
\end{proof}

\begin{corollary} For $\beta = 0$ we have the following special cases.
\begin{itemize}
\item[(i)] $\beta{\mathbb{Y}}=\mathbb{Y}$ is the self-dual graded Young graph.
\item[(ii)] $\varkappa{\mathbb{Y}}$ gives the {\it Pieri deformation} of Young's graph: up edges are as in the usual Young's graph $\mathbb{Y}$ and down edges are given by $(\lambda \to \mu)$ iff $\lambda/\mu$ is a horizontal strip.
\end{itemize}
\end{corollary}

As it was mentioned in \cite{pp}, apparently the most interesting and mysterious type of dual filtered graphs is the so-called {\it M\"obius deformation} that is related to K-theoretic insertion and LR rules. For  $\mathbb{Y}$ it is defined as follows. %We define it more generally as a M\"obius deformation $\mu{\mathbb{Y}}$ of $\mathbb{Y}$. 
%Let $\beta = 1$. (Everything can be done with a general $\beta$ with slight modifications.) 
The defining conditions (i) and (ii) of the M\"obius deformation $\mu\mathbb{Y}$ are the same as for the $\beta$-filtration $\beta{\mathbb{Y}}$ but loops have {\it positive} weight $1$, and {\it down} edges are given by $(\lambda \to \mu)$ iff $\lambda/\mu$ is a {\it rook strip} (i.e., no two boxes lie on the same row or column) with the corresponding weight $w = 1$. % \beta^{|\lambda/\mu|}

Besides new examples of dual filtered Young's graphs, another consequence of our approach is the following result: {M\"obius deformation of Young's lattice is related to the Cauchy deformation and can be obtained from it via a natural transformation}. In particular, this result reveals the presence of a M\"obius deformation for Young's lattice and the transform is in fact related to the M\"obius inversion.

%The following observation shows that one construction gives rise to another relevant filtration or that they are in fact equivalent via this transform.
\begin{lemma}\label{trans}
Suppose $[D, U] = -(1 + D)$. Then $[\widehat{D}, U] = 1 + \widehat{D}$ for $\widehat{D} = -D (1+D)^{-1}$.
\end{lemma}
\begin{proof} See the Appendix.
\end{proof}

\begin{theorem}
The M\"obius deformation  $\mu{\mathbb{Y}}$  is a dual filtered graph satisfying $[D,U] = 1 + D$ and it can be obtained from the Cauchy deformation $\varkappa\mathbb{Y}$ for $\varkappa = \beta = -1$ via the map 
$$
D \longmapsto -D(1 + D)^{-1}. %\qquad \beta \mapsto -\beta.
$$
\end{theorem}
\begin{proof}
Recall that for $\beta = \varkappa = -1$, the down operator of $\varkappa\mathbb{Y}$ (eq. \eqref{opers}) is given by 
\begin{align*}
\overline{D} &= -1 + (1 - \widetilde{d}_1)(1 - \widetilde{d}_2) \cdots  \\
&=-1 + \left(1 -\frac{d_1}{1 + d_1}\right)\left(1 -\frac{d_1}{1 + d_1}\right) \cdots \\
&=-1 + \frac{1}{(1 + d_1)} \frac{1}{(1 + d_2)} \cdots
\end{align*}
Hence,
$$
\widehat{D} = \frac{-\overline{D}}{1 + \overline{D}} = \frac{1}{1 + \overline{D}} - 1 = \cdots (1 + d_2)(1 + d_1) - 1.
$$
Notice that $ \widehat{D} \lambda = (\cdots (1 + d_2)(1 + d_1) - 1) \cdot \lambda$ removes rook strips from $\lambda$ in all possible ways. Therefore, $\widehat{D}$ is a down operator of $\mu\mathbb{Y}$. For $\varkappa\mathbb{Y}$ we have $[\overline{D}, \widetilde{U}] = -(1 + \overline{D})$ and by Lemma~\ref{trans} we have $[\widehat{D}, \widetilde{U}] = 1 + \widehat{D}$.
\end{proof}

\section{Enumerative identities}
%For a dual filtered graph $G$ satisfying the  relation $[D, U] = 1 + q D$ of the {\it $q$-Weyl algebra}, let $f_{G,n}(\lambda/\mu), F_{G,n}(\lambda/\mu)$ be the number of up and respectively down walks of length $n$ from $\mu$ to $\lambda$ in $G$.
Define {\it increasing set-valued tableaux} (ISVT) as an SVT that if after replacing each set by any of its element, the resulting tableau is %{\it standard}, i.e., it 
increasing both in rows and columns. Let $F_{\lambda/\!\!/\mu}(n)$ be the number of ISVT of shape $\lambda/\!\!/\mu$ that contain all numbers from $[n] := \{1,\ldots, n \}$. 

A {\it strict tableaux} (ST) of skew shape $\lambda/\mu$ is a filling of a Young diagram of $\lambda/\mu$ by positive integers so that entries strictly increase in rows from left to right, weakly increase from top to bottom, and each element can appear only on a single column. Let $f_{\lambda/\mu}(n)$ be the number of ST of shape $\lambda/\mu$ that contain all numbers from $[n]$.

An {\it increasing tableaux} (IT) is a filling of a skew diagram by positive integers so that they strictly increase in both rows and columns. Let now $g_{\lambda/\mu}(n)$ be the number of IT of skew shape $\lambda/\mu$ that contain all numbers $[n]$ (some numbers may appear several times).

\begin{theorem} We have
\begin{align*}
\sum_{\lambda} (-1)^{m - |\lambda/\mu|}F_{\lambda/\!\!/\mu}(m) f_{\lambda/\nu}(n) &= \sum_{i} i! \binom{m}{i} \binom{n}{i} \sum_{\kappa} (-1)^{n - i - |\nu/\kappa|} F_{\nu/\!\!/\kappa}(m - i) f_{\mu/ \kappa}(n - i),\\
\sum_{\lambda} F_{\lambda/\!\!/\mu}(m) g_{\lambda/\nu}(n) &= \sum_{i,j} q_{n}(i,j) \binom{m}{i} \binom{n}{j} \sum_{\kappa} F_{\nu/\!\!/\kappa}(m - j) g_{\mu/\kappa}(n - i),
\end{align*}
where
\begin{align}\label{qnij}
q_{n}(i,j) := \sum_{\ell} \binom{i - j + \ell}{\ell} A_{i, n - \ell}
\end{align}
and $A_{i, s}$ is the Eulerian number, i.e., the number of permutations of $(1, \ldots, i)$ with $s$ descents.
\end{theorem}
\begin{proof}
Note that $(-1)^{m - |\lambda/\mu|} F_{\lambda/\!\!/\mu}(m)$ is equivalently the number of {\it signed} {\it up} walks from $\mu$ to $\lambda$ of length $m$ in the $1$-filtration $\beta{\mathbb{Y}}$, i.e., $\beta = 1$, loops have weight $-1$, and a sign of a walk is negative if it uses an odd number of loops, otherwise it is positive. Similarly, $f^{}_{\lambda/\mu}(n)$ is the number of {\it down} walks from $\lambda$ to $\mu$ of length $n$ in $\beta{\mathbb{Y}}$. For example, $F^{}_{(21), (2)}(2) = -3,$ $f_{(211), (1)}(2) = 2$. Enumerator of the down graph of $\mu\mathbb{Y}$ is $g_{\lambda/\mu}$.  Then the formulas are applications of the normal ordering of differential operators $U, D$ given in the lemma below.  
Then the graph $\beta\mathbb{Y}$ with $\beta = 1$ satisfies $[D, U] = 1$ and from Lemma \ref{normord} we obtain:
\begin{align*}
\sum_{\lambda} (-1)^{m - |\lambda/\mu|}F_{\lambda/\!\!/\mu}(m) f_{\lambda/\nu}(n) &= \langle D^n U^m  \mu, \nu \rangle \\
&= \sum_{i} i! \binom{m}{i} \binom{n}{i} \langle U^{m-i} D^{n-i} \mu, \nu \rangle\\
&= \sum_{i} i! \binom{m}{i} \binom{n}{i} \sum_{\kappa} (-1)^{n - i - |\nu/\kappa|} F_{\nu/\!\!/\kappa}(m - i) f_{\mu/ \kappa}(n - i).
\end{align*}
The second formula can be obtained in the same way using the normal ordering for $[D,U] = 1 + D$.
\end{proof}
These formulas are %applications of the normal ordering of differential operators $U, D$ and 
analogous to the formula (see \cite{SS, stadiff, EC2})
$$
\sum_{{|\lambda/\mu| = n \atop |\lambda/\nu| = m}} f_{\lambda/\mu} f_{\lambda/\nu} = \sum_{i \ge 0} i! \binom{m}{i} \binom{n}{i} \sum_{{|\nu/\kappa| = n-i \atop |\mu/\kappa| = m - i}} f_{\nu/\kappa} f_{\mu/\kappa},
$$
where $f_{\lambda/\mu}$ is the number of standard Young tableaux (SYT) of shape $\lambda/\mu$. It generalizes the classical Frobenius identity
$$
\sum_{\lambda\vdash n} f_{\lambda}^2 = n!
$$ 
From the second identity %given in the Theorem is derived from the M\"obius filtration $\mu{\mathbb{Y}}$ and the normal ordering for $[D,U] = 1+ D$, for $\mu = \nu = \varnothing$ 
one can also obtain the formula given in \cite{pp} 
$$
\sum_{\lambda} F_{\lambda}(n) g_{\lambda}(n) = \# \text{ ordered set partitions of } [n].
$$
\begin{lemma}[Normal ordering]\label{normord} The following ordering formulas hold:
\begin{align*}
[D, U] &= 1 \qquad\implies  D^n U^m = \sum_{i} i! \binom{m}{i} \binom{n}{i} U^{m - i} D^{n - i}.\\
[D, U] &= 1 + D \implies D^{n} U^{m} = \sum_{i,j} q_n(i,j) \binom{m}{i} \binom{n}{j} U^{m - i} D^{n - j}. %, \quad q_n(i,j) = \sum_{\ell} \binom{i - j + \ell}{\ell} A_{i, n - \ell},
\end{align*}
%where $A_{i, s}$ is the Eulerian number, i.e., the number of permutations of $(1, \ldots, i)$ with $s$ descents.
%\end{itemize}
\end{lemma}
\begin{proof}
See the Appendix.
\end{proof}

\newpage
\appendix
%\small
\section{Proofs of Lemmas}

\begin{proof}[Proof of Lemma \ref{propcom}]
Non-local identities (i) follow directly from non-local identities in Lemma~\ref{lcom1}. Let us prove the local identities. First, 
{%\small
\begin{align*}
&\qquad\quad \widetilde{u}_{i+1} \widetilde{u}_i (\widetilde{u}_{i+1} + \widetilde{u}_i) = (\widetilde{u}_{i+1} + \widetilde{u}_i) \widetilde{u}_{i+1} \widetilde{u}_i \\
&\iff (u_{i+1} - \beta u_{i+1}d_{i+1}) (u_{i} - \beta u_{i}d_i) (u_{i+1} - \beta u_{i+1}d_{i+1} + u_{i} - \beta u_{i}d_i) \\
&\qquad\qquad\qquad\qquad\qquad\qquad\quad =  (u_{i+1} - \beta u_{i+1}d_{i+1} + u_{i} - \beta u_{i}d_i)(u_{i+1} - \beta u_{i+1} d_{i+1}) (u_{i} - \beta u_{i}d_i) \\
%&\iff u_{i+1} u_i (u_{i+1} + u_i) - \beta u_{i+1} u_i (u_{i+1} d_{i+1} + u_i d_i + d_i u_{i+1} + d_i u_i + d_{i+1} u_{i+1} + d_{i+1} u_i)\\
%&\qquad\quad + \beta^2
\end{align*}
The free coefficients (at $\beta^0$) from both sides are 
$u_{i+1} u_i (u_{i+1} + u_i)$ and $ (u_{i+1} + u_i)u_{i+1} u_i$ that are equal. 
The coefficients at $-\beta$ are equal iff 
\begin{align*}
&\qquad\quad u_{i+1} u_i (\underline{u_{i+1} d_{i+1}} + u_i d_i + d_i u_{i+1} + d_i u_i + d_{i+1} u_{i+1} + \underline{d_{i+1} u_i})\\
&\qquad\qquad\qquad\qquad\qquad=(u_{i+1} d_{i+1} + u_i d_i) u_{i+1}u_i + (u_{i+1} + u_i) u_{i+1}(\underline{d_{i+1} u_i} + u_i d_i) \\
&\iff \underline{u_{i+1} u_i (u_{i+1} + u_i) d_{i+1}} + u_{i+1} u_i (u_i d_i + d_i u_{i+1} + d_i u_i + d_{i+1} u_{i+1})\\
&\qquad\qquad\qquad\qquad\qquad=(u_{i+1} d_{i+1} + u_i d_i) u_{i+1}u_i + \underline{(u_{i+1} + u_i) u_{i+1}u_id_{i+1}} + (u_{i+1} + u_i) u_{i+1} u_i d_i\\
&\iff u_{i+1} u_i (u_i d_i + d_i u_{i+1} + d_i u_i + d_{i+1} u_{i+1})\\
&\qquad\qquad\qquad\qquad\qquad
=(u_{i+1} d_{i+1} + u_i d_i) u_{i+1}u_i + (u_{i+1} + u_i) u_{i+1} u_i d_i\\
&\iff \underline{u_{i+1} u_i (u_i + u_{i+1}) d_i} + u_{i+1} u_i (d_i u_i + d_{i+1} u_{i+1})\\
&\qquad\qquad\qquad\qquad\qquad=(u_{i+1} d_{i+1} + u_i d_i) u_{i+1}u_i + \underline{(u_{i+1} + u_i) u_{i+1} u_i d_i}\\
&\iff u_{i+1} u_i (d_i u_i + d_{i+1} u_{i+1})=(u_{i+1} d_{i+1} + u_i d_i) u_{i+1}u_i \\
&\iff \underline{u_{i+1} u_i d_i u_i} + \underbrace{u_{i+1} u_i d_{i+1} u_{i+1}}_{=u_{i+1} u_i u_i d_i}= \underline{u_{i+1} d_{i+1}u_{i+1}u_i} + u_i d_i u_{i+1}u_i \\
&\iff u_{i+1} u_i u_i d_i = u_i d_i u_{i+1} u_i.
\end{align*}
It is easy to check that the last identity $[u_{i+1} u_i, u_i d_i] = 0$ is always true on the basis elements.

The coefficients at $\beta^2$ are equal iff
\begin{align*}
&\qquad\quad u_{i+1} u_i ((d_i + d_{i+1}) (u_{i+1} d_{i+1} + u_i d_i) + d_{i+1} d_{i}(u_{i+1} + u_i)) \\
&\qquad\qquad\qquad\qquad\qquad= (u_{i+1} + u_i) u_{i+1} u_i d_{i+1} d_i + (u_{i+1} d_{i+1} + u_i d_i) u_{i+1} u_i (d_i + d_{i+1})
\end{align*}
Note that $u_{i+1} d_{i+1} = d_{i+2} u_{i+2}$ commutes with $u_i$ and $d_i$. Let us match the monomials on the l.h.s with the monomials on the r.h.s so that they are equal. We have 

$u_{i+1} u_i d_{i+1} u_i d_i = u_{i+1} d_{i+1} u_i u_i d_i = u_i u_{i+1} d_{i+1} u_i d_i = u_i u_{i+1} u_i d_{i+1} d_i,$

$u_{i+1} u_i d_{i+1} u_{i+1} d_{i+1} = u_{i+1} d_{i+1} u_i u_{i+1} d_{i+1} = u_{i+1} d_{i+1} u_{i+1} d_{i+1} u_i = u_{i+1} d_{i+1} u_{i+1} u_i d_{i+1},$

$
u_{i+1}u_i d_i u_{i+1} d_{i+1} = u_{i+1} u_{i+1} d_{i+1} u_i d_i =  u_{i+1} u_{i+1} u_i d_{i+1} d_i,
$

$
u_{i+1} u_i d_i u_i d_i = u_{i+1} d_{i+1} u_{i+1} u_i d_i,
$

$
u_{i+1} u_i d_{i+1} d_i u_{i +1} = u_{i+1} d_{i+1} u_i d_i u_{i +1} = u_i d_i u_{i+1} d_{i+1} u_{i+1} = u_i d_i u_{i+1} u_i d_i,
$

$
u_{i+1} u_i d_{i+1} d_i u_{i} =u_{i+1} d_{i+1} u_i d_i u_i = u_i d_i u_{i+1} d_{i+1} u_i = u_i d_i u_{i+1} u_i d_{i+1}
$

The coefficients at $-\beta^3$ are equal iff
\begin{align*}
&\qquad\quad u_{i+1}d_{i+1} u_i d_i (u_{i+1}d_{i+1} + u_i d_i) =  (u_{i+1}d_{i+1} + u_i d_i) u_{i+1}d_{i+1} u_i d_i
\end{align*}
and since the elements $\{ u_i d_i \}$ commute, $[u_i d_i, u_{i+1} d_{i+1}] = [u_i d_i, d_{i+2} u_{i+2}] = 0$, the identity is true.

Let us prove now the local identity for $\widetilde{d}$. The identity 
\begin{align*}
\widetilde{d}_{i} \widetilde{d}_{i+1} (\widetilde{d}_{i} + \widetilde{d}_{i+1}) = (\widetilde{d}_{i} + \widetilde{d}_{i+1})\widetilde{d}_{i} \widetilde{d}_{i+1} 
\end{align*}
is equivalent to 
\begin{align}\label{sumb}
\sum_{k,\ell,m \ge 1} \beta^{k + \ell + m - 3} d_{i}^{k}  d_{i+1}^{\ell} (d_{i}^{m} + d_{i+1}^m) = \sum_{k',\ell',m' \ge 1} \beta^{k' + \ell' + m' - 3} (d_{i}^{m'} + d_{i+1}^{m'})  d_{i}^{k'}  d_{i+1}^{\ell'}.
\end{align}
Since $[d_i d_{i+1}, d_i] = [d_{i} d_{i+1}, d_{i+1}] =0$ %$[d_{i}^{k} d_{i+1}^{\ell}, d_{i}^{m}] = [d_{i}^{k} d_{i+1}^{\ell}, d_{i+1}^{m}] = 0$ 
and for all $k, \ell, m \ge 1$ we have 
$$
d_{i}^{k} d_{i+1}^{\ell}  d_{i}^{m} = d_i^{k - 1} d_i d_{i+1} d_{i+1}^{\ell - 1} d_{i}^m %= d_i^{k - 1}  d_{i+1}^{\ell - 1} d_{i}^m d_i d_{i+1} 
= d_i^{k - 1}  d_{i+1}^{\ell - 1} d_{i}^{m+1} d_{i+1} = \cdots = d_{i}^{k - t} d_{i+1}^{\ell - t} d_{i}^{m+t} d_{i+1}^t,
$$
where $ t = \min(k, \ell)$. If  $t = \ell \le k$ then $d_{i}^{k} d_{i+1}^{\ell}  d_{i}^{m} = d_{i}^{m}d_i^{k} d_{i+1}^\ell$ and if $t = k < \ell$ then $d_{i}^{k} d_{i+1}^{\ell}  d_{i}^{m} = d_{i+1}^{\ell - k} d_{i}^{m + k} d_{i+1}^k$. Observe that $(d_i d_{i+1})^s = d_{i}^s d_{i+1}^s$ and then
$$
d_{i}^{k} d_{i+1}^{\ell}  d_{i+1}^{m} = d_{i}^{k - t} (d_{i} d_{i+1})^t d_{i+1}^{\ell - t} d_{i+1}^m = d_{i}^{k - t} d_{i+1}^{m} d_{i}^t d_{i+1}^{\ell}. %= d_{i}^k d_{i+1}^{m + \ell} = d_{i}^{k-1} d_i d_{i+1} d_{i+1}^{m + \ell - 1} = d_{i}^{k-1} d_{i+1}^{m + \ell -1 } d_{i} d_{i+1} = \cdots = d_{i}^{k - t} d_{i+1}^{m + \ell - t} d_i^t d_{i+1}^t
$$
Again, if $t = k \le \ell$ then $d_{i}^{k} d_{i+1}^{\ell}  d_{i+1}^{m} = d_{i+1}^m d_{i}^k d_{i+1}^{\ell}$ and if $t = \ell < k$ then $d_{i}^{k} d_{i+1}^{\ell}  d_{i+1}^{m} = d_{i}^{k - \ell} d_{i+1}^{m} d_{i}^\ell d_{i+1}^{\ell} = d_{i}^{k - \ell} d_{i}^{\ell} d_{i+1}^{m + \ell}$. 
Therefore, every element of the l.h.s. of \eqref{sumb} can be matched with the elements of the r.h.s. as follows
\begin{align*}
d_{i}^{k}  d_{i+1}^{\ell} d_{i}^{m} =  d_{i}^{m}d_i^{k} d_{i+1}^\ell, \qquad \text{ if } \ell \le k\\
d_{i}^{k}  d_{i+1}^{\ell} d_{i}^{m} =   d_{i+1}^{\ell - k} d_{i}^{m + k} d_{i+1}^k, \qquad \text{ if } \ell > k\\
d_{i}^{k} d_{i+1}^{\ell}  d_{i+1}^{m} =   d_{i+1}^m d_{i}^k d_{i+1}^{\ell}, \qquad \text{ if } k \le \ell\\
d_{i}^{k} d_{i+1}^{\ell}  d_{i+1}^{m} =   d_{i}^{k - \ell} d_{i}^{\ell} d_{i+1}^{m + \ell}, \qquad \text{ if } k > \ell
\end{align*}
Note that the degree $k + \ell + m - 3$ is preserved and it is easy to check that we have defined a $\beta$ degree preserving bijection between the elements of l.h.s. and r.h.s of \eqref{sumb}. 

Let us finally verify the conjugate relations. The relation $[\widetilde{u}_i, \widetilde{d}_j] = 0$ for $|i - j| \ge 2$ is an easy consequence of a non-local commutativity. Now $\widetilde{d}_1 \widetilde{u}_1 = (1- \beta d_1)^{-1} d_1 u_1 (1 - \beta d_1) = 1$ and 
$$[\widetilde{u}_{i+1}, \widetilde{d}_i] = [u_{i+1} - \beta u_{i+1} d_{i+1}, d_i (1 - \beta d_i)^{-1}] = [u_{i+1} - \beta d_{i+2} u_{i+2}, d_i (1 - \beta d_i)^{-1}] = 0.$$
}
\end{proof}

\begin{proof}[Proof of Lemma \ref{locud}]
Expanding the identity 
$$
(1 - xy \widetilde{u}_i \widetilde{d}_i)^{-1} (1 + x\widetilde{u}_i)(1+y \widetilde{d}_{i+1}) = 
(1 - xy \widetilde{d}_{i+1} \widetilde{u}_{i+1})^{-1} (1 + y\widetilde{d}_{i+1})(1+x \widetilde{u}_{i}).  \eqno{(*)}
$$
we need to show that
$$
\sum_{k} (xy)^{k} (\widetilde{u}_i \widetilde{d}_i)^k (1 + x\widetilde{u}_i + y \widetilde{d}_{i+1} + xy \widetilde{u}_i \widetilde{d}_{i+1}) = \sum_{k} (xy)^k (\widetilde{d}_{i+1} \widetilde{u}_{i+1})^k (1 + x\widetilde{u}_i + y \widetilde{d}_{i+1} + xy \widetilde{d}_{i+1} \widetilde{u}_{i})
$$
or, equivalently, that the following identities hold (the first one for $k \ge 1$ and others for $k \ge 0$)
\begin{align}
(\widetilde{u}_i \widetilde{d}_i)^k + (\widetilde{u}_i \widetilde{d}_i)^{k-1} \widetilde{u}_i \widetilde{d}_{i+1} &= (\widetilde{d}_{i+1} \widetilde{u}_{i+1})^k + (\widetilde{d}_{i+1} \widetilde{u}_{i+1})^{k-1} \widetilde{d}_{i+1} \widetilde{u}_{i} \label{ud1} \\
(\widetilde{u}_i \widetilde{d}_i)^k \widetilde{d}_{i+1} &= (\widetilde{d}_{i+1} \widetilde{u}_{i+1})^k \widetilde{d}_{i+1} \label{ud2}\\
(\widetilde{u}_i \widetilde{d}_i)^k \widetilde{u}_{i} &= (\widetilde{d}_{i+1} \widetilde{u}_{i+1})^k \widetilde{u}_{i} \label{ud3}
\end{align}
Recall that $\widetilde{u}_i =u_i(1 - \beta d_i) = u_i - \beta u_i d_{i}$, $\widetilde{d}_{i} = d_{i}(1 - \beta d_{i})^{-1} = \sum_{\ell \ge 1} \beta^{\ell - 1} d_{i}^{\ell}$, and hence $\widetilde{u}_i \widetilde{d}_i = u_i d_i$.  Recall also that $u_i d_i = d_{i+1} u_{i+1}$ and $[u_i, d_j] = 0$ for $i \ne j$.

Let us show \eqref{ud1}. For $k = 1$ we have 
{\small
\begin{align}
&\qquad\quad \widetilde{u}_i \widetilde{d}_i + \widetilde{u}_i \widetilde{d}_{i+1} = \widetilde{d}_{i+1} \widetilde{u}_{i+1} + \widetilde{d}_{i+1} \widetilde{u}_{i}  \label{udu}\\
&\iff u_i d_i + (u_i - \beta u_i d_i)\sum_{\ell \ge 1} \beta^{\ell - 1} d_{i+1}^{\ell} = \sum_{\ell \ge 1} \beta^{\ell - 1} d_{i+1}^{\ell} (u_{i+1} - \beta u_{i+1}d_{i+1} + u_i - \beta u_i d_{i}) \nonumber\\
&\iff u_i d_i + \sum_{\ell \ge 1} \underline{\beta^{\ell - 1} u_i d_{i+1}^{\ell}} - \beta^{\ell} u_i d_i d_{i+1}^{\ell} = \sum_{\ell \ge 1} \beta^{\ell - 1} d_{i+1}^{\ell} u_{i+1} - \beta^{\ell} d_{i+1}^{\ell} u_{i+1} d_{i+1} + \underline{\beta^{\ell - 1} d_{i+1}^{\ell} u_i} - \beta^{\ell} d_{i+1}^{\ell} u_i d_i \nonumber\\
&\iff 
u_{i}d_{i} - \sum_{\ell \ge 1} \beta^{\ell} u_{i} d_{i} d_{i+1}^{\ell} = \sum_{\ell \ge 1} \underline{\beta^{\ell - 1} d_{i+1}^{\ell} u_{i+1}} - \beta^{\ell} d_{i+1}^{\ell} u_{i+1} d_{i+1}  - \underline{\beta^{\ell} d_{i+1}^{\ell} d_{i+1} u_{i+1}} \nonumber\\
&\iff 
u_{i}d_{i} - \sum_{\ell \ge 1} \beta^{\ell} u_{i} d_{i} d_{i+1}^{\ell} = d_{i+1} u_{i+1} - \sum_{\ell \ge 1}\beta^{\ell} d_{i+1}^{\ell} u_{i+1} d_{i+1} \nonumber\\
&\iff 
\sum_{\ell \ge 1} \beta^{\ell} u_{i} d_{i} d_{i+1}^{\ell} = \sum_{\ell \ge 1}\beta^{\ell} d_{i+1}^{\ell} u_{i+1} d_{i+1} \nonumber
\end{align}
} 
Notice that from $[d_i d_{i+1}, d_{i+1}] = 0$ we have $$u_i d_i d_{i+1}^{\ell} = u_i d_{i+1}^{\ell - 1} d_i d_{i+1} = d_{i+1}^{\ell - 1} u_i d_i d_{i+1} =  d_{i+1}^{\ell - 1} d_{i+1} u_{i+1} d_{i+1} = d_{i+1}^{\ell} d_{i+1} d_{i+1}$$
meaning that the preceding identities are indeed true. Next, note that $(\widetilde{u}_i \widetilde{d}_i)^k = (u_i d_i)^{k} = u_i d_i = \widetilde{u}_i\widetilde{d}_i$. For $k = 2$ we need to show that 
{\small
\begin{align*}
&\qquad\ (\widetilde{u}_i \widetilde{d}_i)^2 + \widetilde{u}_i \widetilde{d}_i \widetilde{u}_i \widetilde{d}_{i+1} = (\widetilde{d}_{i+1} \widetilde{u}_{i+1})^2 + \widetilde{d}_{i+1} \widetilde{u}_{i+1} \widetilde{d}_{i+1} \widetilde{u}_{i}\\
&\iff
\widetilde{u}_i \widetilde{d}_i + \widetilde{u}_i \widetilde{d}_i \widetilde{u}_i \widetilde{d}_{i+1} = (\widetilde{d}_{i+1} \widetilde{u}_{i+1})^2 + \widetilde{d}_{i+1} \widetilde{u}_{i+1} \widetilde{d}_{i+1} \widetilde{u}_{i}\\
&\iff
u_i d_i + u_i d_i \widetilde{u}_i \widetilde{d}_{i+1} = \widetilde{d}_{i+1} u_{i+1} d_{i+1} \widetilde{u}_{i+1} + \widetilde{d}_{i+1} u_{i+1} d_{i+1} \widetilde{u}_i\\
&\iff 
u_i d_i + u_i d_i(u_i - \beta u_i d_i) \sum_{\ell \ge 1} \beta^{\ell - 1} d_{i+1}^{\ell} = \sum_{\ell \ge 1} \beta^{\ell - 1} d_{i+1}^{\ell} u_{i+1} d_{i+1}(u_{i+1} - \beta u_{i+1} d_{i+1} + u_i - \beta u_i d_i)\\
&\iff 
u_i d_i +  \sum_{\ell \ge 1} \beta^{\ell - 1} u_i d_i u_i d_{i+1}^{\ell} - \beta^{\ell}  u_i d_i  d_{i+1}^{\ell} = \sum_{\ell \ge 1} \underline{\beta^{\ell - 1} d_{i+1}^{\ell} u_{i+1} d_{i+1}u_{i+1}} - \beta^{\ell} d_{i+1}^{\ell} u_{i+1} d_{i+1} u_{i+1} d_{i+1}\\ 
&\qquad \qquad\qquad\qquad\qquad\qquad\qquad\qquad\qquad\qquad\qquad+ \beta^{\ell - 1} d_{i+1}^{\ell} u_{i+1} d_{i+1} u_i - \underline{\beta^{\ell} d_{i+1}^{\ell} u_{i+1} d_{i+1} u_i d_i}
%&\iff
%u_i d_i + u_i d_i {u}_i (1 - \beta d_i) {d}_{i+1}(1 - \beta d_{i+1})^{-1} = (1 - \beta d_{i+1})^{-1} {d}_{i+1} u_{i+1} d_{i+1} {u}_{i+1} (1 - \beta d_{i+1}) \\
%&\qquad\qquad\qquad\qquad\qquad\qquad\qquad\qquad\qquad\qquad\qquad + (1- \beta d_{i+1})^{-1}{d}_{i+1} u_{i+1} d_{i+1} {u}_i(1 - \beta d_i)\\
%&\iff
%u_i d_i + u_i d_i {u}_i (1 - \beta d_i) {d}_{i+1}(1 - \beta d_{i+1})^{-1} = (1 - \beta d_{i+1})^{-1} u_i d_i (1 - \beta d_{i+1}) \\
%&\qquad\qquad\qquad\qquad\qquad\qquad\qquad\qquad\qquad\qquad\qquad + (1- \beta d_{i+1})^{-1}u_i d_i d_{i+1} {u}_i(1 - \beta d_i)
\end{align*}
} 
Notice that %$d_{i+1} u_{i+1} d_{i+1} d_{i+1} u_{i+1} = 
$u_i d_i d_{i+1} u_i d_i = d_{i+1} u_i d_i$ (removal of a box in the column $i+1$ does not change the property if the box in the column $i$ was removable). Hence for $\ell \ge 2$
{\small
$$d_{i+1}^{\ell} u_{i+1} d_{i+1}u_{i+1} = d_{i+1}^{\ell-1} u_i d_i u_i d_i = d_{i+1}^{\ell-1} u_i d_i = d_{i+1}^{\ell - 2} d_{i+1} u_i d_i = d_{i+1}^{\ell - 2} u_i d_i d_{i+1} u_i d_i =  d_{i+1}^{\ell - 2} d_{i+1} u_{i+1} d_{i+1} u_i d_i $$
}
and then the last sum is equivalent to the following simplified form
{\small
\begin{align*}
&\iff 
u_i d_i +  \sum_{\ell \ge 1} \beta^{\ell - 1} u_i d_i u_i d_{i+1}^{\ell} - \beta^{\ell}  u_i d_i  d_{i+1}^{\ell} = \underbrace{d_{i+1} u_{i+1} d_{i+1} u_{i+1}}_{=u_id_i} \\
&\qquad+ \sum_{\ell \ge 1} - \beta^{\ell} \underbrace{d_{i+1}^{\ell} u_{i+1} d_{i+1} u_{i+1} d_{i+1}}_{=d_{i+1}^{\ell-1} u_i d_i d_{i+1} = u_i d_{i+1}^{\ell-1} d_i d_{i+1} = u_i d_i d_{i+1}^{\ell}} + \beta^{\ell - 1} \underbrace{d_{i+1}^{\ell} u_{i+1} d_{i+1} u_i}_{=d_{i+1}^{\ell - 1} u_i d_i d_{i+1} u_i = u_i d_{i+1}^{\ell - 1} d_i d_{i+1} u_i = u_i d_i d_{i+1}^{\ell} u_i = u_i d_i u_i d_{i+1}^{\ell} }
\end{align*}
}
which is true. Now for $k > 2$ first recall that 
$(\widetilde{u}_i \widetilde{d}_i)^k = \widetilde{u}_i \widetilde{d}_i$ and hence $$(\widetilde{d}_{i+1} \widetilde{u}_{i+1})^k =\widetilde{d}_{i+1} (\widetilde{u}_{i+1} \widetilde{d}_{i+1})^{k-1} \widetilde{u}_{i+1} = (\widetilde{d}_{i+1} \widetilde{u}_{i+1})^2$$
and so we obtain
\begin{align*}
(\widetilde{u}_i \widetilde{d}_i)^k + (\widetilde{u}_i \widetilde{d}_i)^{k-1} \widetilde{u}_i \widetilde{d}_{i+1} &= \widetilde{u}_i \widetilde{d}_i + \widetilde{u}_i \widetilde{d}_i \widetilde{u}_i \widetilde{d}_{i+1}\\
&=(\widetilde{d}_{i+1} \widetilde{u}_{i+1})^2 + \widetilde{d}_{i+1} \widetilde{u}_{i+1} \widetilde{d}_{i+1} \widetilde{u}_{i} \qquad\qquad \text{(using $k = 2$)}\\
&=(\widetilde{d}_{i+1} \widetilde{u}_{i+1})^k + \widetilde{d}_{i+1} (\widetilde{u}_{i+1} \widetilde{d}_{i+1})^{k-1} \widetilde{u}_{i} \\
&=(\widetilde{d}_{i+1} \widetilde{u}_{i+1})^k + (\widetilde{d}_{i+1} \widetilde{u}_{i+1})^{k-1} \widetilde{d}_{i+1} \widetilde{u}_{i} \label{ud1}
\end{align*}
which gives \eqref{ud1}.

Let us now prove \eqref{ud2}. Again, for $k = 1$ we need to prove that
\begin{align*}
&\qquad\quad \widetilde{u}_i \widetilde{d}_i \widetilde{d}_{i+1} = \widetilde{d}_{i+1} \widetilde{u}_{i+1} \widetilde{d}_{i+1} \\
&\iff u_i d_i d_{i+1} (1 - \beta d_{i+1})^{-1} = (1 - \beta d_{i+1})^{-1} d_{i+1} u_{i+1} d_{i+1}\\
&\iff u_i d_i d_{i+1} (1 - \beta d_{i+1})^{-1} = (1 - \beta d_{i+1})^{-1} u_i d_i d_{i+1}\\
&\iff u_i d_i d_{i+1} (1 - \beta d_{i+1})^{-1} = u_i (1 - \beta d_{i+1})^{-1} d_i d_{i+1}
\end{align*}
which is true since $[d_i d_{i+1}, d_{i+1}^{\ell}] = 0$. For $k \ge 2$ we have
$$
(\widetilde{u}_i \widetilde{d}_i)^k \widetilde{d}_{i+1} =  \widetilde{u}_i \widetilde{d}_i \widetilde{d}_{i+1} = \widetilde{d}_{i+1} \widetilde{u}_{i+1} \widetilde{d}_{i+1} = \widetilde{d}_{i+1} (\widetilde{u}_{i+1} \widetilde{d}_{i+1})^k = (\widetilde{d}_{i+1} \widetilde{u}_{i+1})^k \widetilde{d}_{i+1}
$$
which gives \eqref{ud2}.

Finally, let us show \eqref{ud3}. For $k = 1$ we need to show that
\begin{align*}
&\qquad\quad\widetilde{u}_{i} \widetilde{d}_{i} \widetilde{u}_{i} = \widetilde{d}_{i+1} \widetilde{u}_{i+1} \widetilde{u}_{i}\\
&\iff u_i d_i u_i (1 - \beta d_i) = (1 - \beta d_{i+1})^{-1} \underbrace{d_{i+1} u_{i+1}}_{=u_id_i} (1 - \beta d_{i+1}) u_i (1 - \beta d_i)\\
&\iff u_i d_i u_i (1 - \beta d_i) = u_i (1 - \beta d_{i+1})^{-1} d_i (1 - \beta d_{i+1}) u_i (1 - \beta d_i)\\
&\iff u_i d_i u_i (1 - \beta d_i) = u_i (1 - \beta d_{i+1})^{-1} \underbrace{d_i u_i}_{=u_{i-1}d_{i-1}} (1 - \beta d_{i+1}) (1 - \beta d_i)\\
&\iff u_i d_i u_i (1 - \beta d_i) = u_i d_i u_i (1 - \beta d_{i+1})^{-1} (1 - \beta d_{i+1}) (1 - \beta d_i)\\
&\iff u_i d_i u_i (1 - \beta d_i) = u_i d_i u_i (1 - \beta d_i).
\end{align*}
For $k = 2$ we need to show that 
\begin{align*}
&\qquad\quad(\widetilde{u}_{i} \widetilde{d}_{i})^2 \widetilde{u}_{i} = (\widetilde{d}_{i+1} \widetilde{u}_{i+1})^2 \widetilde{u}_{i}\\
&\iff u_i d_i u_i (1 - \beta d_i) = (1 - \beta d_{i+1})^{-1} \underbrace{d_{i+1} u_{i+1} d_{i+1} u_{i+1}}_{=u_i d_i} (1 - \beta d_{i+1}) u_i (1 - \beta d_i) 
\end{align*}
and then it follows the same way as the previous chain for $k = 1$. Now for $k > 2$ we have
$$
(\widetilde{u}_{i} \widetilde{d}_{i})^k \widetilde{u}_{i} = \widetilde{u}_{i} \widetilde{d}_{i} \widetilde{u}_{i} = (\widetilde{d}_{i+1} \widetilde{u}_{i+1})^2 \widetilde{u}_{i} = \widetilde{d}_{i+1} (\widetilde{u}_{i+1} \widetilde{d}_{i+1})^{k-1} \widetilde{u}_{i+1} \widetilde{u}_{i} = (\widetilde{d}_{i+1} \widetilde{u}_{i+1})^k \widetilde{u}_{i}
$$
which gives \eqref{ud2}.
\end{proof}

\begin{proof}[Proof of Lemma \ref{trans}]
Note that $\widehat{D} = (1 + D)^{-1} - 1$. We have 
\begin{align*}
&DU - UD = -(1 + D)\\ 
\implies &(1+D)U - U(1+D) = -(1 + D) \\
\implies &(1+D)U(1+D)^{-1} - U = -1\\
\implies &U(1+D)^{-1} - (1+D)^{-1}U = -(1+D)^{-1}\\
\implies &U((1+D)^{-1}-1) - ((1+D)^{-1}-1)U = -(1+D)^{-1} + 1 - 1\\
\implies &U\widehat{D} - \widehat{D}U = -\widehat{D} - 1\\
\implies &[\widehat{D},U] = 1 + \widehat{D}.
\end{align*}
\end{proof}

\begin{proof}[Proof of Lemma \ref{normord}]
The part (i) is standard and well-known.
To prove (ii) first observe that $(1 + D)U - U(1 + D) = 1+ D$ and then by induction $(1 + D)^n U = (U + n) (1 + D)^n$. Now we have
\begin{align*}
D^n U^m &= ((1 + D) - 1)^n U^m\\
		& = \sum_{k} (-1)^{n - k} \binom{n}{k} (1 + D)^k u^m\\
		& = \sum_{k} (-1)^{n - k} \binom{n}{k}  (U + k)^m (1 + D)^k\\
		& = \sum_{k} (-1)^{n - k} \binom{n}{k} \sum_i \binom{m}{i} k^{m - i} U^{i} \sum_{j} \binom{k}{j} D^j\\
		& = \sum_{i,j} U^i D^j \sum_{k} (-1)^{n - k} \binom{n}{k} \binom{m}{i} \binom{k}{j} k^{m - i}\\
		& = \sum_{i,j} \binom{m}{i} \binom{n}{j} U^i D^j \sum_{k} (-1)^{n - k} \binom{n-j}{k -j} k^{m-i}\\
		& = \sum_{i,j} \binom{m}{i} \binom{n}{j} U^{m-i} D^{n-j} \sum_{k} (-1)^{n - k} \binom{j}{n -k} k^i.
\end{align*}
Note that $q_n(i,j) =  \sum_{k} (-1)^{n - k} \binom{j}{n -k} k^i$ is the coefficient at $[z^n]$ of the series 
$$
\sum_{\ell} (-1)^{\ell} \binom{j}{\ell} z^\ell \sum_{m} z^m m^i = (1 - z)^j \frac{A_i(z)}{(1 - z)^{i+1}} = (1 - z)^{i - j + 1} A_i(z),
$$
where $A_i(z) = \sum_{s = 0}^{i - 1} A_{i,s} z^{s}$ is the Eulerian polynomial.  
And thus $q_n(i,j) = \sum_{\ell} \binom{i - j + \ell}{\ell} A_{i, n - \ell}$ as needed.
\end{proof}

%\end{appendix}

\end{document}